\documentclass[12pt]{article}

\usepackage{amsfonts,amsmath, amssymb}
\usepackage[utf8]{inputenc}
\usepackage{amsmath}
\usepackage{amssymb}
\usepackage{amscd}
\usepackage{amsfonts}
\usepackage{amsthm}
\usepackage{mathrsfs}
\usepackage[usenames, dvipsnames]{color}
\usepackage[matrix,arrow,curve]{xy}

\usepackage{xcolor}

\usepackage{enumerate}
\usepackage{epsf,epsfig,amsfonts,a4wide}
\usepackage{amsfonts, mathdots}
\usepackage{amsmath,amssymb,amsthm}
\usepackage{url}
\usepackage{amsmath,amssymb,amsthm}

\usepackage{amssymb}
\usepackage{amsthm}
\usepackage{epsf,epsfig,amsfonts,graphicx}
\usepackage{a4wide}
\newtheorem{theor}{Theorem}[section]
\newtheorem{lemma}{Lemma}[section]
\newtheorem{prop}{Proposition}[section]
\newtheorem{corr}{Corollary}[section]

\newtheorem{example}{Example}[section]

\theoremstyle{remark}

\newcommand{\be}{\begin{equation}}
\newcommand{\ee}{\end{equation}}

\newcommand{\sfm}{\mathsf{m}}
\newcommand{\sfl}{\mathsf{l}}
\newcommand{\gl}{\mathfrak{gl}}
\newcommand{\ddet}{\operatorname{det}}

\newcommand{\ddd}{\mathrm{d}}

\newcommand{\pd}[2]{\frac{\partial#1}{\partial#2}}

\newcommand{\tr}{\operatorname{tr}}

\newcommand{\weg}[1]{}

\title{Nijenhuis geometry II: Left-symmetric algebras and linearization problem for Nijenhuis operators}
\author{Andrey Yu. Konyaev\footnote{Faculty of Mechanics and Mathematics, Moscow State University, 119992, Moscow Russia
 \ \ \quad {\tt  maodzund@yandex.ru}}}  
\date{}

\begin{document}

\maketitle

\section{Introduction}

The Nijenhuis torsion \cite{nijenhuis} of an operator field $R$ on a manifold is a tensor defined on a pair of vector fields $v, w$ as:
\begin{equation}\label{torsion}
    \mathcal N_R(v, w) = R[Rv, w] + R[v, Rw] - R^2[v, w] - [Rv, Rw].
\end{equation}
Here $[\, , \,]$ stands for the standard commutator of vector fields. The vanishing of the Nijenhuis torsion \eqref{torsion} in local coordinates $x^1, \dots, x^n$ is equivalent to the system of partial differential equations on the components $R^{\alpha}_i$ of $R$: for $i \neq j$ and $1 \leq i, j, k \leq n$
\begin{equation}\label{torsion2}
    \pd{R^{\alpha}_i}{x^j} R_{\alpha}^k - \pd{R^{\alpha}_j}{x^i} R_{\alpha}^k -\pd{R^k_i}{x^{\alpha}} R^{\alpha}_j + \pd{R^k_j}{x^{\alpha}} R^{\alpha}_i = 0.
\end{equation}
Operator field $R$ with vanishing Nijenhuis torsion is called \textbf{Nijenhuis operator}. Nijenhuis operators play important role in finite and infinite dimensional integrable systems \cite{kossman}, \cite{mok2}, almost complex structures \cite{newlander} and projectively equivalent metrics \cite{proj}.

The point $\mathrm{p}$ is \textbf{algebraically generic} if Segre characteristic of $R$ at $\mathrm{p}$ is the same as Segre characteristic of $R$ for all points in sufficiently small neighbourhood of $\mathrm{p}$. If point is not algebraically generic, then it is \textbf{singular}. The point $\mathrm{p}$ is \textbf{a singular point of scalar type} if $R = \lambda \operatorname{Id}$ for some constant $\lambda$ at this point. The following Splitting Theorem \cite{bmk} play an important role in study of the local geometry of the Nijenhuis operators:
\begin{theor} Assume that the spectrum of a Nijenhuis operator $R$ at a
point ${\mathrm{p}}$ consists of $k$ real (distinct) eigenvalues $\lambda_1,\dots, \lambda_k$  with multiplicities $\sfm_1,\dots, \sfm_k$  respectively and $s$ pairs of complex (non-real) conjugate eigenvalues $\mu_1,\bar \mu_1, \dots ,\mu_s,\bar \mu_s$ of multiplicities $\sfl_1,\dots, \sfl_s$.  Then in a neighborhood of ${\mathrm{p}}$ there exists a local coordinate system 
$$
x_1=(x_{1}^1\dots  x_{1}^{\sfm_1}),  \dots  , x_k=(x_{k}^1 \dots  x_{k}^{\sfm_k}),\quad 
u_1=(u_{1}^1 \dots  u_{1}^{2\sfl_1}),  \dots  , u_s=(u_{s}^1 \dots  u_{s}^{2\sfl_s}),
$$
in which $R$ takes the following block-diagonal form
$$
R = 
\begin{pmatrix} 
Q_1(x_1)&  & & & & \\
 &\ddots & & &  &\\
 & & Q_k(x_k)& & & \\
 & & & Q^{\mathbb C}_1 (u_1)&  & \\
 & & & & \ddots  & \\
 & & & & & Q^{\mathbb C}_s (u_s)
\end{pmatrix}
$$
where each block depends on its own group of variables and is a Nijenhuis operator w.r.t. these variables.   
\end{theor}

This theorem reduces the study of the Nijenhuis operators around any point to the study of the Nijenhuis operators around the point, where all eigenvalues of the given operator coincide. The singular points of scalar type are the most singular of such points.

In the present paper we follow the Nijenhuis geometry programme, formulated in \cite{bmk}, and study the local properties of the Nijenhuis operators in the neighbourhood of the points of scalar type. Recall, that all of the classical results about Nijenhuis operators \cite{nijenhuis}, \cite{haantjes}, \cite{thompson}, \cite{newlander} were obtained for algebraically generic points, often with extra regularity conditions.

Before formulating the results of the paper we need a couple of definitions. Let $\mathfrak a$ be an algebra of dimension $n$ over $\mathbb R$ with the multiplication $\star$. The associator $\mathcal A$ is a trillinear operation on $\mathfrak a$, defined on arbitrary triple $\xi, \eta, \zeta \in \mathfrak a$ as $\mathcal A(\xi, \eta, \zeta) = (\xi\star\eta)\star\zeta - \xi\star(\eta\star\zeta)$ for $\xi, \eta, \zeta \in \mathfrak a$. Associator identically vanish if and only if $\mathfrak a$ is associative. An algebra $\mathfrak a$ is called \textbf{left-symmetric or LSA} if its associator is symmetric in first two terms, that is 
\begin{equation}\label{lsa}
\mathcal A (\xi, \eta , \zeta) = \mathcal A (\eta, \xi, \zeta), \quad \forall \xi, \eta, \zeta \in \mathfrak a.    
\end{equation}
The main property of these algebras is, that commutator $[\xi, \eta] = \xi\star\eta - \eta\star\xi$ defines a Lie algebra structure on $\mathfrak a$. We call this Lie algebra \textbf{the associated Lie algebra}. The other names for these algebras, that appear in the literature, are pre-Lie algebras and Vinberg-Kozul algebras. These algebras were introduced by Vinberg \cite{vinberg} in his study of the homogeneous cones. Later they appeared in different frameworks of geometry, integrable systems and quantum mechanics (see \cite{burde} for an overview on the subject). In particular, the Novikov algebras, that play an important role in the theory of Poisson brackets of the hydrodynamic type are left-symmetric. 

The left-adjoint action $L_{\xi}$ and the right-adjoint $R_{\xi}$ of $\mathfrak a$ on itself are defined as usual: $L_{\xi} \eta = \xi \star \eta = R_{\eta} \xi$. In terms of the left-adjoint action the condition \eqref{lsa} can be written as
\begin{equation}\label{left}
L_{\xi} L_{\eta} - L_{\eta} L_{\xi} = L_{[\xi, \eta]}.
\end{equation}

Assume now, that $\mathfrak a$ is an arbitrary finite dimensional algebra over $\mathbb R$ of dimension $n$. Every such algebra has a natural structure of a smooth $n-$dimensional affine manifold. We denote manifold same as algebra itself. Consider point $\eta$ of this manifold. Tangent space to $\eta$ is naturally identified with $\mathfrak a$ itself. Thus, one can define a tensor field $R$ of type $(1, 1)$ as follows: $R$ acts on element $\xi$ from tangent space $T_{\eta} \mathfrak a$ by the right-adjoint action, that is $R_{\eta} \xi = \xi\star\eta$.

Now fix a basis $\eta_i$ in $\mathfrak a$ and denote by $a_{ij}^k$ the structure constants of $\mathfrak a$. The basis defines natural coordinate system on $\mathfrak a$ and denote the corresponding coordinates as $x^i$, that is $\eta = x^i \eta_i$. The components of $R_{\eta}$ are written as $(R_{\eta})^k_i = a^k_{ij} x^j$. In particular, the entries are homogeneous linear functions of $x^i$. 

Thus, the right-adjoint action of $\mathfrak a$ on itself induces the operator field with homogeneous linear components. We call operator field $R_{\eta}$ \textbf{the right-adjoint operator of $\mathfrak a$}. In a similar fashion one constructs \textbf{left-adjoint operator of $\mathfrak a$}. We denote it same as left action as $L_{\eta}$.

Moreover, assume one has an operator field $R_{\eta}$ on a real affine space $\mathfrak a$ with given coordinates $x^i$ and its entries are homogeneous linear functions (we call such operator fields \textbf{linear operator fields}). Then $\mathfrak a$ has a natural structure of algebra over $\mathbb R$ and its structure constants are $a^k_{ij} = \pd{R^k_i}{x^j}$. In this construction $R$ becomes the right-adjoint operator of $\mathfrak a$. This yields a natural bijection between real algebras and linear operator fields on real affine spaces. 

We call linear operator field on affine space with vanishing Nijenhuis torsion \textbf{linear Nijenhuis operator}. The following Proposition establishes the relation between linear Nijenhuis operators and left-symmetric algebras. For the sake of simplicity we omit $\eta$ in $R_{\eta}, L_{\eta}$ and just write $R, L$ respectively. 

\begin{prop}\label{main1}(\cite{winterhalder}, Theorem 3.1, Remark 3.1)
Let $\mathfrak a$ be an algebra over $\mathbb R$ of dimension $n$. The following conditions are equivalent:
\begin{enumerate}
    \item $\mathfrak a$ is a left-symmetric algebra
    \item The right-adjoint operator of $\mathfrak a$ is a Nijenhuis operator
\end{enumerate}
\end{prop}

Note, also, that left-adjoint operator $L$ is not Nijenhuis in general.

\begin{example}\label{dir}
Consider a linear Nijenhuis operator 
\begin{equation*}
    R = \left(\begin{array}{cccc}
    x^1 & 0 & \dots & 0 \\
    0 & x^2 & \dots & 0 \\
     & & \ddots & \\
     0 & 0 & \dots & x^n \\
    \end{array}\right).
\end{equation*}
The corresponding left-symmetric algebra $\mathfrak a$ in corresponding basis $\eta_1, \dots, \eta_n$ has structure constants:
\begin{equation*}
    \begin{aligned}
        a^k_{ij} = \begin{cases}
        1, \quad i = j = k, \\
        0, \quad \text{\rm{otherwise}}
        \end{cases}
    \end{aligned}
\end{equation*}
The Lie algebra, associated to $\mathfrak a$, is abelian and $\mathfrak a$ is a direct sum of one-dimensional left-symmetric algebras with operation $\star$ defined as $\eta \star \eta = \eta$. 
\end{example}

Now we are proceed to the results of the paper. First, we classify all real left-symmetric algebras in dimension two. Until now there has been only a partial ($\mathfrak b$-series in our terminology) classification by Burde (Proposition, 4.1, \cite{burde2}) of the left symmetric algebras over $\mathbb C$. In our classification we adopt notation, similar to the one in \cite{burde2}. For $\mathfrak b_2, \mathfrak b^+_4, \mathfrak b^-_4$ we have slightly different basis and $\mathfrak b^+_4, \mathfrak b^-_4$ have the same complexification, denoted in \cite{burde2} as $\mathfrak b_4$. 

\begin{theor}\label{class}
Up to isomorphism there are two continuous families and 10 exceptional two dimensional real left-symmetric algebras. The complete list of normal forms is presented in the Table 1 and Table 2 below. The tables contain 
\begin{itemize}
    \item All non-zero structure relations for a given basis $\eta_1, \eta_2$
    \item The right-adjoint operator of $\mathfrak a$ in coordinates $x, y$, associated with basis $\eta_1, \eta_2$. We denote it as $R$
    \item The left-adjoint operator of $\mathfrak a$ in coordinates $x, y$, associated with basis $\eta_1, \eta_2$. We denote it as $L$.
\end{itemize}
Recall, that in dimension two up to the isomorphism there are only two Lie algebras. The letter $\mathfrak b$ stands for algebras with non-abelian associated Lie algebra and $\mathfrak c$ for the algebras with abelian associated Lie algebra.
\begin{equation*}
\begin{aligned}
& Table \, 1: \text{LSA in dimension two with non-commutative associated Lie algebra} \\
& \begin{tabular}{|l|l|l|l|}
\hline
Name & Structure relations & $L$ & $R$ \\
\hline
$\mathfrak b_{1, \alpha}$ & $\begin{array}{ll} & \eta_2\star\eta_1 = \eta_1,\\ & \eta_2\star\eta_2 = \alpha \eta_2\\ \end{array}$ & $\left(\begin{array}{cc} y & 0 \\ 0 & \alpha y\\ \end{array}\right)$ & $\left(\begin{array}{cc} 0 & x \\ 0 & \alpha y\\ \end{array}\right)$ \\
\hline
$\begin{array}{ll} & \mathfrak b_{2, \beta} \\ & \beta \neq 0\\ \end{array}$ & $\begin{array}{ll} & \eta_1\star\eta_2 = \eta_1, \\& \eta_2\star\eta_1 = \big(1 - \frac{1}{\beta}\big) \eta_1 \\ & \eta_2\star\eta_2 = \eta_2 \end{array}$ &
$\left(\begin{array}{cc} \big(1 - \frac{1}{\beta}\big) y & x \\ 0 & y\\ \end{array}\right)$ & $\left(\begin{array}{cc} y & \big(1 - \frac{1}{\beta}\big) x \\ 0 & y\\ \end{array}\right)$\\
\hline
$\mathfrak b_3$ & $\begin{array}{ll} & \eta_2\star\eta_1 = \eta_1, \\ &  \eta_2\star\eta_2 = \eta_1 + \eta_2\\ \end{array}$ & $\left(\begin{array}{cc} y & y \\ 0 & y\\ \end{array}\right)$ & $\left(\begin{array}{cc} 0 & x + y \\ 0 & y\\ \end{array}\right)$ \\
\hline
$\mathfrak b_4^{+}$ & $\begin{array}{ll} & \eta_1\star\eta_1 = \eta_2, \\ &  \eta_2\star\eta_1 = - \eta_1\\ & \eta_2\star\eta_2 = - 2 \eta_2 \end{array}$ & $\left(\begin{array}{cc} - y & 0 \\ x & - 2y\\ \end{array}\right)$ & $\left(\begin{array}{cc} 0 & -x \\ x & - 2y\\ \end{array}\right)$ \\
\hline
$\mathfrak b_4^{-}$ & $\begin{array}{ll} & \eta_1\star\eta_1 = - \eta_2, \\ &  \eta_2\star\eta_1 = - \eta_1\\ & \eta_2\star\eta_2 = - 2 \eta_2 \end{array}$ & $\left(\begin{array}{cc} - y & 0 \\ - x & - 2y\\ \end{array}\right)$ & $\left(\begin{array}{cc} 0 & - x \\ - x & -2y\\ \end{array}\right)$ \\
\hline
$\mathfrak b_5$ & $\begin{array}{ll} & \eta_1\star\eta_2 = \eta_1, \\ &  \eta_2\star\eta_2 = \eta_1 + \eta_2\\ \end{array}$ & $\left(\begin{array}{cc} 0 & x + y \\ 0 & y\\ \end{array}\right)$ & $\left(\begin{array}{cc} y & y \\ 0 & y\\ \end{array}\right)$ \\
\hline
\end{tabular}
\end{aligned}
\end{equation*}
\begin{equation*}
\begin{aligned}
& Table \, 2: \text{LSA in dimension two with commutative associated Lie algebra} \\
& \begin{tabular}{|l|l|l|}
\hline
Name & Structure relations & $L = R$\\ 
\hline
$\mathfrak c_1$ &  &  $\left(\begin{array}{cc} 0 & 0 \\ 0 & 0\\ \end{array}\right)$ \\ 
\hline
$\mathfrak c_2$ & $\eta_2\star\eta_2 = \eta_2$ & $\left(\begin{array}{cc} 0 & 0 \\ 0 & y\\ \end{array}\right)$ \\ 
\hline
$\mathfrak c_3$ & $\eta_2\star\eta_2 = \eta_1 $ & $\left(\begin{array}{cc} 0 & y \\ 0 & 0\\ \end{array}\right)$ \\ 
\hline
$\mathfrak c_4$ & $\begin{aligned}& \eta_2\star\eta_2 = \eta_2\\ & \eta_2\star\eta_1 = \eta_1 \\& \eta_1\star\eta_2 = \eta_1 \\ \end{aligned}$ & $\left(\begin{array}{cc} y & x \\ 0 & y\\ \end{array}\right)$ \\ 
\hline
$\mathfrak c_5^{+}$ & $\begin{aligned}& \eta_2\star\eta_2 = \eta_2\\ & \eta_2\star\eta_1 = \eta_1 \\& \eta_1\star\eta_2 = \eta_1 \\ & \eta_1\star\eta_1 = \eta_2\\ \end{aligned}$ & $\left(\begin{array}{cc} y & x \\ x & y\\ \end{array}\right)$ \\ \hline
$\mathfrak c_5^{-}$ & $\begin{aligned}& \eta_2\star\eta_2 = \eta_2\\ & \eta_2\star\eta_1 = \eta_1 \\& \eta_1\star\eta_2 = \eta_1 \\ &\eta_1\star\eta_1 = - \eta_2\\ \end{aligned}$ & $\left(\begin{array}{cc} y & x \\ - x & y\\ \end{array}\right)$ \\ \hline
\end{tabular}
\end{aligned}
\end{equation*}
\end{theor}

Fix singular point $\mathrm{p}$ of scalar type of Nijenhuis operator $R$ on manifold $M^n$ and coordinates $x^1, \dots, x^n$ with coordinate origin is at $\mathrm{p}$. Consider a pair of vectors $v_{\mathrm{p}}, w_{\mathrm{p}} \in T_{\mathrm{p}} M^n$. Denote by $v, w$ arbitrary smooth continuations of $v_{\mathrm{p}}, w_{\mathrm{p}}$ on the neighborhood of ${\mathrm{p}}$. Define the following operation
\begin{equation}\label{isotropy}
v_{\mathrm{p}} \star w_{\mathrm{p}} = \big ([Rv, w] - R [v, w] \big)\vert_{\mathrm{p}}.
\end{equation}
The following proposition is the second result of the paper:

\begin{prop}\label{main2}
Let $R$ be a Nijenhuis operator and ${\mathrm{p}}$ singular point of scalar type. Then 
\begin{enumerate}
    \item The definition of $\star$ does not depend on continuations $v, w$ and defines an algebra structure on $T_{\mathrm{p}} M^n$ 
    \item The corresponding algebra is left-symmetric
    \item In natural basis of $T_{\mathrm{p}} M^n$, associated with coordinates $x^1, \dots, x^n$, the structure constants of this algebra are $a^k_{ij} = \pd{R^k_i}{x^j}\vert_{\mathrm{p}}$.
\end{enumerate}
\end{prop}

The corresponding algebra $\mathfrak a$ is called \textbf{the isotropy algebra at ${\mathrm{p}}$}. Note, that if $R$ is Nijenhuis, then $R - \lambda \operatorname{Id}$ for constant $\lambda$ is also Nijenhuis. Moreover, if ${\mathrm{p}}$ is a singular point of scalar type for $R$, then ${\mathrm{p}}$ is a singular point of scalar type for $R - \lambda \operatorname{Id}$ and in given coordinates the isotropy algebras coincide. We will often assume, that $R$ vanishes at ${\mathrm{p}}$.

Tangent space $T_{\mathrm{p}} M^n$ has a natural structure of affine space and $\mathfrak a$ defines by Proposition \ref{main1} natural Nijenhuis operator on it --- the right-adjoint operator. Thus, we have two Nijenhuis operators: one is on tangent space $T_{\mathrm{p}} M^n$ and the other is on the manifold, defined in a neighbourhood of $\mathrm{p}$. 

W.l.o.g. assume that $R$ vanishes at ${\mathrm{p}}$ and consider the Taylor expansion of $R$ at ${\mathrm{p}}$: $R_1 + R_{k + 1} + \dots$. The entries of $R_i$ are homogeneous polynomials of degree $i$. By Proposition \ref{main2} the first term of the expansion $R_1$ in given coordinates is:
$$
(R_1)^k_i = \pd{R^k_i}{x^j}\vert_{\mathrm{p}} x^j.
$$
The term $R_1$ is not an operator field on $M^n$ --- it behaves wrong under the coordinate change. But in every coordinate system it is a result of simple lifting of the right-adjoint operator from $T_{\mathrm{p}}$ onto the neighbourhood of ${\mathrm{p}}$ simply by replacing the affine coordinates on $T_{\mathrm{p}} M^n$ with the coordinates on $M^n$. 

So the natural question arises: is there a coordinate system, where this lifting yields the entire $R$? Or, in other words, is there a coordinate change, in which $R$ coincides $R_1$ on the entire neighbourhood of ${\mathrm{p}}$? This is \textbf{the linearization problem for Nijenhuis operators}. Following the terminology of Weinstein \cite{weinstein}, \cite{weinstein2} we call left-symmetric algebra $\mathfrak a$ \textbf{non-degenerate} if following property holds: if isotropy algebra of Nijenhuis operator $R$ at a singular point of scalar type ${\mathrm{p}}$ is isomorphic to $\mathfrak a$, then there exists a linearizing coordinate change. 

In \cite{bmk} it is proved that the left-symmetric algebra, described in Example \ref{dir}, is non-degenerate in both formal and analytic category. In this work we give the complete classification of real two-dimensional left-symmetric algebras in terms of non-degeneracy in the smooth category and incomplete classification in the real analytic category (this is the third result of our paper).

Introduce the following subsets of $\mathbb R$
$$
\Sigma_0 = \{0\}, \quad \Sigma_1 = \{r \vert \, r \in \mathbb N, r \geq 3\}, \quad \Sigma_2 = \{\alpha \vert \, \alpha \in \mathbb R, \alpha < 0\}, \quad \Sigma_3 = \{\frac{1}{m} \vert \, m \in \mathbb N, r \geq 2\}
$$
and denote $\Sigma_{\mathrm{sm}} = \Sigma_0 \cup \Sigma_1 \cup \Sigma_2 \cup \Sigma_3$. The following Theorem holds:

\begin{theor}\label{themain1}
The following table provides the complete classification of two-dimensional left-symmetric algebras in terms of non-degeneracy in the smooth category:
\begin{equation*}
\begin{aligned}
& Table \, 3: \text{Classification of LSA in dimension two in smooth category} \\
& \begin{tabular}{|c|c|}
\hline
    Degenerate LSA & Non-degenerate LSA  \\
     \hline
     $\begin{array}{c}\mathfrak c_1, \mathfrak c_2, \mathfrak c_3, \mathfrak c_4, \\ \mathfrak b_5, \mathfrak b_{2, \beta}\\
     \mathfrak b_{1, \alpha} \, \text{\rm{for}} \,  \alpha \in \Sigma_{\mathrm{sm}}
     \end{array}$ & $\begin{array}{c} \mathfrak b^+_4, \mathfrak b^-_4, \mathfrak c^+_5, \mathfrak c^-_5 \\ \mathfrak b_3, \mathfrak b_{1, \alpha} \, \text{\rm{for}} \, \alpha \notin \Sigma_{\mathrm{sm}}\end{array}$ \\
     \hline
\end{tabular}\\
\end{aligned}
\end{equation*}
\end{theor}

Let $[q_0, q_1, q_2, ...]$ be a decomposition of an irrational $\alpha$ into the continuous fraction. If the series
\begin{equation*}
B(x) = \sum \limits_{i = 0}^{\infty}\frac{\operatorname{log}q_{i + 1}}{q_i}
\end{equation*}
converges, then $\alpha$ is a \textbf{Brjuno number}. We denote by $\Omega$ \textbf{the set of negative Brjuno numbers} and $\Sigma_{\mathrm{u}}$ is the set of negative irrational numbers, that are not Brjuno numbers. Note, that the Lebesque measure of $\Sigma_{\mathrm{u}}$ is zero.

Define $\widehat{\Sigma}_2 = \{- \frac{p}{q}\ \vert \, p, q \in \mathbb N\}$ and consider $\Sigma_{\mathrm{an}} = \Sigma_0 \cup \Sigma_1 \cup \widehat{\Sigma}_2 \cup \Sigma_3$. The following Theorem holds.

\begin{theor}\label{themain2}
The following table provides the (incomplete) classification of two-dimensional left-symmetric algebras in terms of non-degeneracy in the analytic category:
\begin{equation*}
\begin{aligned}
& Table \, 4: \text{Classification of LSA in dimension two in analytic category} \\
& \begin{tabular}{|c|c|c|}
\hline
    Degenerate LSA & Non-degenerate LSA & Unknown \\
     \hline
     $\begin{array}{c}
     \mathfrak c_1, \mathfrak c_2, \mathfrak c_3, \mathfrak c_4, \\ 
     \mathfrak b_5, \mathfrak b_{2, \beta}
     \mathfrak b_{1, \alpha}\, for \, \alpha \in \Sigma_{\mathrm{an}}
     \end{array}$ & $\begin{array}{c} 
     \mathfrak b^+_4, \mathfrak b^-_4, \mathfrak c^+_5, \mathfrak c^-_5 \\
     \mathfrak b_3, \mathfrak b_{1, \alpha} \, for \, \alpha \notin \Sigma_{\mathrm{an}} \cup \Sigma_{\mathrm{u}}
     \end{array}$ & $\begin{array}{c}
          \mathfrak b_{1, \alpha} \, for \, \alpha \in \Sigma_{\mathrm{u}}
     \end{array}$ \\
     \hline
\end{tabular} \\
\end{aligned}
\end{equation*}
\end{theor}

The set $\Sigma_{\mathrm{u}}$ has an interesting story, directly related to the linearization problem of $R$ in dimension two. Consider analytic vector field $v$ with critical point ${\mathrm{p}}$ and denote $\lambda_1, \lambda_2$ the eigenvalues of the linearization operator at ${\mathrm{p}}$. We are interested in linearization problem for such vector field (see Appendix A for details). If $\frac{\lambda_1}{\lambda_2}$ is negative rational number or $\frac{\lambda_1}{\lambda_2} = r, \frac{1}{r}$ for $2 \leq r \in \mathbb N$, then there are no formal linearization and, thus, no analytic one. This is the classical result by variety of authors, including Poincare himself (see Theorem \ref{ilyash}). If $\frac{\lambda_1}{\lambda_2}> 0$ and $\frac{\lambda_1}{\lambda_2} \neq r, \frac{1}{r}$ for $2 \leq r \in \mathbb N$, then there exists an analytic linearizing coordinate change (see Theorem \ref{poincare}). For $\alpha \in \Sigma$ the linearization problem was solved by Brjuno (see Theorem \ref{brjn}). The question remained open only for $\frac{\lambda_1}{\lambda_2} \in \Sigma_{\mathrm{u}}$. 

In the late 80s---early 90s this gap was closed by Yoccoz \cite{yoccoz1}, \cite{yoccoz2}, who showed, that for every $\alpha \in \Sigma_{\mathrm{u}}$ there exists such an analytic vector field with $\frac{\lambda_1}{\lambda_2} = \alpha$, that among all formal linearizing coordinate changes, there are no converging ones. Thus, these are no analytic linearizing coordinate change for such vector fields. 

In our case the only "troublesome" algebra is $\mathfrak b_{1, \alpha}$. At some point the linearization problem for this algebra is reduced to the linearization (see Lemma \ref{prenormal}) of certain vector field $v$ in a special form $v = (f(x, y), y)^T$ (we write vector fields as transposed rows). The linearizing coordinate change must also be in a special form $\bar x = g(x, y), \bar y = y$. We use the mentioned above results to either provide an example of Nijenhuis operator without linearizing coordinate change or prove the existence of such linearization. Unfortunately, for $\alpha \in \Sigma_u$ the "bad" analytic vector field, Yoccoz constructs in his work, is obtained from the complex one, thus it is not in the triangular form we want. Thus, we were not able to apply his result and the our classification in analytic category is incomplete. It seems, that the linearization problem in triangular form was not considered by ODE specialists before.

The work is organized as follows: Section \ref{local_geometry} provides some analytical results about Nijenhuis operators in dimension two. All the necessary definitions and results from the theory of ODE, concering the linearization problem and existence of the first integrals around critical points, are presented in Appendix A. Appendix B contains the proof of so called Morse lemma depending on parameters and its Corollary. We included it to keep the work self-sufficient. The proof of Theorem \ref{class} and Proposition \ref{main2} are given in sections \ref{proof1} and \ref{proof2}. The proof of Theorems \ref{themain1} and \ref{themain2} is given in section \ref{proof3}.

Author would like to thank Ilya Schurov, Yuri Kudryashov, Alexei Bolsinov for help and fruitful discussions. He would also like to thank referees for valuable comments and suggestions. The work was supported by Russian Science Foundation (project No. 17-11-01303).


\section{Local geometry of Nijenhuis operators in dimension 2}\label{local_geometry}

Fix coordinates $x, y$ and consider Nijenhuis operator 
\begin{equation*}
    R = \left(\begin{array}{cc}
    R^1_1 & R^1_2 \\
    R^2_1 & R^2_2 \\
    \end{array}\right).
\end{equation*}

The following proposition holds:
\begin{prop}\label{local1}
In dimension two the following conditions are equivalent:\\
1) Operator $R$ is Nijenhuis;\\
2) In local coordinates $x, y$:
\begin{equation}\label{id1}
\ddd \ddet R = \ddd \tr R \, \left(\begin{array}{cc}
R^2_2 &  - R^1_2 \\
- R^2_1 & R^1_1
\end{array}\right) \,,
\end{equation}
where $\operatorname{det} R, \operatorname{tr} R$ are the determinant and trace of $R$ respectively. 
\end{prop}
{\it Proof.} From \eqref{torsion2} for $i = 1, j = 2, k = 1$ and $x^1 = x, x^2 = y$ we get 
\begin{equation*}
\begin{aligned}
    0 & = \pd{R^1_1}{y} R^1_1 + \pd{R^2_1}{y} R^1_2 - \underline{\pd{R^1_2}{x} R_1^1} - \pd{R^2_2}{x} R^1_2 - \\
    & - \pd{R^1_1}{x}R^1_2 - \pd{R^1_1}{y}R^2_2 + \underline{\pd{R^1_2}{x}R^1_1} + \pd{R^1_2}{y}R^2_1 
\end{aligned}
\end{equation*}
Underlined terms cancel. Adding $\pd{R^2_2}{y}R^1_1$ to r.h.s. and subtracting it we get
\begin{equation*}
    \begin{aligned}
    0 & = \pd{}{y} \Big( - R^1_1 R^2_2 + R^2_1 R^1_2\Big) - R^1_2 \Big( \pd{R^1_1}{x} + \pd{R^2_2}{x}\Big) + R^1_1 \Big( \pd{R^1_1}{y} + \pd{R^2_2}{y}\Big).
    \end{aligned}
\end{equation*}
In a similar way for $i = 1, j = 2, k = 2$ formula \eqref{torsion2} yields
\begin{equation*}
    \begin{aligned}
    0 & = \pd{}{x} \Big( - R^1_1 R^2_2 + R^2_1 R^1_2\Big) + R^2_2 \Big( \pd{R^1_1}{x} + \pd{R^2_2}{x}\Big) - R^2_1 \Big( \pd{R^1_1}{y} + \pd{R^2_2}{y}\Big).
    \end{aligned}
\end{equation*}
As $R^1_1 R^2_2 - R^2_1 R^1_2 = \ddet R$ and $R^1_1 + R^2_2 = \tr R$ these equations can be rewritten in form \eqref{id1}. The Lemma is proved. $\blacksquare$\\

\begin{corr}\label{useful}
Assume, that $f(x, y)$ and $g(y)$ are arbitrary (smooth or analytic) functions of two and one variables respectively. It follows from Theorem \ref{local1} that
\begin{equation*}
R = \left(
\begin{array}{cc}
g(y) & f(x, y)\\
0 & g(y)
\end{array}\right)
\end{equation*}
is Nijenhuis operator.
\end{corr}

\begin{corr}\label{useful2}
Assume, that $f(x, y)$ and $g(y)$ are arbitrary (smooth or analytic) functions of two and one variables respectively. It follows from Theorem \ref{local1} that
\begin{equation*}
R = \left(
\begin{array}{cc}
     0 & f(x, y)  \\
     0 & g(y) 
\end{array}\right),
\end{equation*}
is Nijenhuis operator.
\end{corr}

Consider Nijenhuis $R$ with condition $\ddd \operatorname{tr} R \neq 0$ at point ${\mathrm{p}}$. Fix $\alpha \neq 0$ and choose coordinates $x, y$ around ${\mathrm{p}}$ such that $\operatorname{tr} R = R^1_1 + R^2_2 = \alpha y$. Denote $\operatorname{det} R$ as $g (x, y)$. Proposition \eqref{id1} yields
\begin{equation}\label{full}
    \begin{aligned}
        R^1_1 &R^2_2 - R^1_2 R^2_1 = g, \\
        R^2_1 & = - \frac{1}{\alpha} \pd{g}{x}, \quad R^1_1 = \frac{1}{\alpha} \pd{g}{y}, \\
        R^1_1 & + R^2_2 = \alpha y. \\
    \end{aligned}
\end{equation}

We treat this system as a system with functional parameter $g$. We get the following corollary.

\begin{corr}\label{local2}
Every solution of system \eqref{full} is written in the form:
\begin{equation}\label{solution}
    R = \left(\begin{array}{cc}
    \frac{1}{\alpha} \pd{g}{y} & R^1_2 \\
    - \frac{1}{\alpha} \pd{g}{x} & y - \frac{1}{\alpha} \pd{g}{y} \\
    \end{array}\right),
\end{equation}
where $R^1_2$ satisfies the following (implicit) condition 
\begin{equation}\label{condition}
\frac{1}{\alpha}\pd{g}{y} \Big (y - \frac{1}{\alpha}\pd{g}{y}\Big) + \frac{1}{\alpha}\pd{g}{x} R^1_2 - g = 0.
\end{equation}
\end{corr}

Condition \eqref{condition} does not contain any derivatives of $R^1_2$. It can be solved for arbitrary function $g$, but the component $R^1_2$ may not even be continuous. 

\begin{example}\label{bath}
Consider the special case: $\alpha = 1$ and $g$ depends only on $y$. In this case $\pd{g}{y} = g'$. Condition \eqref{condition} yields
$$
\Big(y - g'\Big) g' - g = 0
$$
Differentiating it in $y$ we get
$$
g''(y - 2 g') = 0.
$$
This implies, that $g$ is a polynomial in $y$ of at most of second degree. Using the method of undetermined coefficients we get two solutions: $g = \frac{y^2}{4}$ and $g = \beta y - \beta^2$ for arbitrary constant $\beta$.
\end{example}

The following proposition further explores the analytical properties of Nijenhuis operators, that follow from formula \eqref{solution}.

\begin{prop}\label{normal}
Consider a real plane with coordinates $x, y$ and analytic (smooth) Nijenhuis operator $R$ with singular point of scalar type at the coordinate origin with $\lambda = 0$. Assume, that $\ddd \operatorname{tr} R \neq 0$ at the origin and determinant $g = \operatorname{det} R$ satisfies the following conditions at the origin:
$$
g(0, 0) = \pd{g}{x}(0, 0) = \pd{g}{y}(0, 0) = \frac{\partial^2 g}{\partial x^2} (0, 0) = \frac{\partial^2 g}{\partial y^2} (0, 0) = \frac{\partial^2 g}{\partial x \partial y} (0, 0) = 0.
$$
In other words, $g$ has no constant, linear and quadratic part in its Taylor expansion. Define constant $\gamma$ as $\pd{R^1_2}{x} (0, 0) = \gamma$. We get
\begin{enumerate}
    \item In analytic category assume that 
    $$
\gamma \neq 0, \quad \gamma \neq - \frac{p}{q} \text{ with } p, q \in \mathbb N, \quad \gamma \neq \frac{1}{m} \text{ with } 3 \leq m \in \mathbb N,
$$
\item In smooth category assume, that
$$
\gamma > 0, \quad \gamma \neq \frac{1}{m} \text{ with } 3 \leq m \in \mathbb N.
$$
\end{enumerate}
Under these assumptions there exists such an analytic (smooth) coordinate change, that in new coordinates $\tilde x, \tilde y$ Nijenhuis operator $R$ takes form:
\begin{equation}\label{normal_normal}
    R = \left ( \begin{array}{cc}
         0 & f(\tilde x, \tilde y)  \\
         0 & \tilde y
    \end{array}\right),
\end{equation}
with $f(0, 0) = 0$ and $\pd{f}{\tilde x} (0, 0) = \gamma$.
\end{prop}
{\it Proof. }We start with lemma.

\begin{lemma}\label{analytic}
Consider ODE
\begin{equation}\label{op}
    r(x)\cdot k'(x) - k(x) = 0
\end{equation}
in a neighbourhood of the point $x = 0$. Assume, that $r(0) = 0$ and $r'(0) = \beta \neq 0$. Then
\begin{enumerate}
    \item If $\frac{1}{\beta} \notin \mathbb N$ then the only analytic (smooth) solution of \eqref{op} is $k(x) \equiv 0$
    \item If $\frac{1}{\beta} \in \mathbb N$ then every analytic (smooth) solution in sufficiently small neighbourhood of $x = 0$ can be written as $k(x) = c x^{\frac{1}{\beta}} F(x)$, where $F(x)$ is analytic (smooth) function, defined only by $r(x)$, $F(0) \neq 0$ and $c$ is arbitrary constant.
\end{enumerate}
\end{lemma}
{\it Proof. } Consider $k(x)$ to be an analytic (smooth) solution of \eqref{op} on the neighbourhood of $x = 0$. W.l.o.g. we assume that this neighbourhood is $x \in [-\epsilon, \epsilon]$ for some $\epsilon > 0$.

As $r'(0) \neq 0$ then in sufficiently small neighbourhood of $x = 0$ the only zero of $r(x)$ is $x = 0$. Thus, w.l.o.g. we assume that $r(x) \neq 0$ for $x \in [-\epsilon, 0) \cup (0, \epsilon]$. Again, as $r'(0) \neq 0$, then by Morse lemma there exists an analytic (smooth) function $\bar r(x)$ such that $r(x) = \beta x + x^2 \bar r(x)$. For $x \in [-\epsilon, 0) \cup (0, \epsilon]$ the equation \eqref{op} can be rewritten as
$$
\frac{k'}{k} = \frac{1}{r(x)} = \frac{1}{\beta x + x^2 \bar r(x)} = \frac{1}{\beta x} - \frac{1}{\beta}\frac{\bar r}{\beta + x \bar r}. 
$$
Assume, that $x \in (0, \epsilon]$. Integrating both sides of the equation and taking the exponent, we get
\begin{equation}\label{form}
k(x) = c x^{\frac{1}{\beta}} \underline{\exp \Big( \int \limits^{\epsilon}_x \frac{1}{\beta}\frac{\bar r(t)}{\beta + t \bar r(t)} \ddd t \Big)} = c x^{\frac{1}{\beta}} F(x).    
\end{equation}
Here $F(x)$ denotes the underlined integral. As $r(x) \neq 0$ for $x \in [-\epsilon, 0) \cup (0, \epsilon]$, we have that $\beta + x \bar r(x) \neq 0$ for $x \in [-\epsilon, \epsilon]$. Thus, the integral, that defines $F(x)$ converges for all $x$ from the domain and $F(x)$ is analytic (smooth) function with property $F(x) \neq 0$. In particular, $F(0) \neq 0$.

Assume now that the first statement of Lemma \ref{analytic} is wrong. That is $\frac{1}{\beta} \notin \mathbb N$ and $c \neq 0$ for some analytic (smooth) solution $k(x)$. Pick the smallest $m \in \mathbb Z^+$ such, that $m > \frac{1}{\beta}$ and consider $m-$th derivative of $k(x)$. We get that $\lim_{x \to 0+} k^{(m)}(x) = \infty$. This contradiction implies, that the first statement of the Lemma is proved.

Now consider $\frac{1}{\beta} \in \mathbb N$. The formula \eqref{form} defines analytic (smooth) function on the entire $x \in [\epsilon, \epsilon]$ and this function coincides with the solution on $x \in (0, \epsilon]$. By uniqueness theorem for ODE we get that $k(x)$ coincides with formula \eqref{form} on the entire domain of $x$ and the second statement of Lemma \eqref{analytic} is proved. $\blacksquare$ 

\begin{lemma}\label{analytic2}
In the condition of the Proposition \ref{normal} there exists an analytic (smooth) function $\mu(x, y)$, such that $\mu(0, 0) = 0, \pd{\mu}{x} (0, 0) = 0, \pd{\mu}{y} (0, 0) = 1$ and
\begin{equation}\label{astra}
    \operatorname{det} \big( R - \mu \operatorname{Id}\big) \equiv 0.
\end{equation}
In other words $\mu$ is eigenvalue of $R$ at each point.
\end{lemma}
{\it Proof. } Function $\operatorname{det} R = g(x, y)$ satisfies the condition of Corollary \eqref{morse2}. Thus, for $\alpha = 1$ there exist analytic (smooth) functions $h(x, y), k(x)$ such, that 
\begin{equation}\label{rt3}
h(0, 0) = \pd{h}{x}(0, 0) = 0, \quad \Big( \pd{h}{y}(0, 0)\Big)^2 = \frac{1}{4} \neq 0, \quad k(0) = k'(0) = k''(0) = 0
\end{equation}
and
\begin{equation}\label{rt2}
    g(x, y) = \frac{y^2}{4} - h^2(x, y) + k(x).
\end{equation}
As $\pd{h}{y} (0, 0) \neq 0$, then by the Implicit Function Theorem in sufficiently small neighbourhood of the origin there exists an analytic (smooth) curve $s(x)$ such, that 
$$
h(x, s(x)) = 0, \quad s(0) = 0, \quad s'(0) = - \pd{h}{x}(0, 0) \big( \pd{h}{y}(0, 0)\big)^{-1} = 0.
$$
The last property follows from \eqref{rt3}.

Substituting \eqref{rt2} into the condition \eqref{condition} for $\alpha = 1$ yields 
$$
\Big( \frac{y}{2} - 2 h \pd{h}{y}\Big)\Big( \frac{y}{2} + 2 h \pd{h}{y}\Big) + \Big(k' - 2 h \pd{h}{x}\Big) R^1_2 - \frac{y^2}{4} + h^2 - k = 0. 
$$
Putting $y = s(x)$ and renaming $R^1_2(x, s(x)) = r(x)$ we get
$$
r(x)k'(x) - k(x) = 0.
$$
At the same time $r'(0) = \pd{R^1_2}{x} (0, 0) + \pd{R^1_2}{y}(0, 0) s'(0) = \gamma$. 

The conditions of Proposition \ref{normal} imply, that $\gamma$ is either $1, \frac{1}{2}$ or $\frac{1}{\gamma} \notin \mathbb N$. If $\frac{1}{\gamma} \notin \mathbb N$ then, by the first statement of Lemma \ref{analytic} $k(x) \equiv 0$. For $\gamma = 1, \frac{1}{2}$ recall, that by properties \eqref{rt3} $k'(0) = k''(0) = 0$. Thus, constant $c$ for $k(x)$ from second statement of the Lemma \ref{analytic} must be zero. Thus, for all $\gamma$ function $k \equiv 0$.

Finally, consider a pair of analytic (smooth) functions $\mu_1 (x, y) = \frac{y}{2} + h(x, y)$ and $\mu_2 (x, y) = \frac{y}{2} - h(x, y)$. As $\mu_1 + \mu_2 = y$ and $\mu_1 \mu_2 = \frac{y^2}{4} - h^2 = g$, then both functions are roots of characteristic polynomial of $R$. Thus, both functions satisfy condition \eqref{astra}. From conditions \eqref{rt3} we know, that $\pd{h}{y} (0, 0)$ is either $\frac{1}{2}$ or $- \frac{1}{2}$. For $\frac{1}{2}$ take $\mu = \mu_1$ and for $- \frac{1}{2}$ take $\mu = \mu_2$. We have that $\pd{\mu}{y} (0, 0) = 1$. By the same conditions \eqref{rt3} $\pd{\mu}{x} (0, 0) = \pm \pd{h}{x} (0, 0) = 0$. The Lemma is proved. $\blacksquare$

The condition \eqref{astra} implies, that $\mu(x, y)$ is an eigenvalue of $R$ at every point. For Nijenhuis operator $R$ such function $\mu$ satisfies the following condition (\cite{bmk}, Proposition 2.3, Formula 9):
$$
R^* \ddd \mu = \mu \ddd \mu.
$$
We take $\mu$, that exists by Lemma \ref{analytic2}, and consider analytic (smooth) coordinate change $\bar x = x, \bar y = \mu(x, y)$. In new coordinates
\begin{equation}\label{rt4}
R = \left( \begin{array}{cc}
     R^1_1 & R^1_2  \\
     0 & \bar y 
\end{array}\right).
\end{equation}
The trace takes the form $\bar y + R^1_1$ and determinant is $g = yR^1_1$. From conditions of Proposition \ref{normal} on $g$ we have:
$$
\pd{g}{\bar y} (0, 0) = R^1_1 (0, 0), \quad \frac{\partial^2 g}{\partial \bar y^2} (0, 0) = 2 \pd{R^1_1}{\bar y} (0, 0), \quad \frac{\partial^2 g}{\partial \bar x \partial \bar y}(0, 0) = \pd{R^1_1}{\bar x} (0, 0).
$$
Also, note, that the Jacobian of the coordinate change $\bar x(x, y), \bar y (x, y)$ at the coordinate origin is $\operatorname{Id}$. Thus, the coordinate change does not change the linear part of Taylor expansion of $R$ at the origin and, thus, $\pd{R^1_2}{x} (0, 0) = \gamma$. 

From Proposition \ref{local1} the vanishing of the Nijenhuis torsion of $R$ in the form \eqref{rt4} yields the single equation
\begin{equation}\label{int}
R^1_2 \pd{R^1_1}{\bar x} + \Big(\bar y - R^1_1 \Big) \pd{R^1_1}{\bar y} = 0.    
\end{equation}
In coordinates $\bar x, \bar y$ consider vector field $v = \big(R^1_2, \bar y - R^1_1\big)^T$. 

We have that the coordinate origin is a critical point of $v$. The eigenvalues of linearization operator at this point are $1, \gamma$. Equation \eqref{int} implies, that $R^1_1$ is a first integral of vector field $v$. By Corollary \ref{nodes} under our assumption on $\gamma$ there are no non-constant analytic (smooth) first integral around the origin for $v$. Thus, $R^1_1 \equiv 0$. The Proposition is proved. $\blacksquare$

\begin{example}
Lemma \ref{analytic2} implies, that the eigenvalue of the Nijenhuis operator under the conditions of Proposition \ref{normal} is analytic (smooth) function. In general, this is not the case. Consider the right-adjoint operator of $\mathfrak b_4^-$: 
\begin{equation}\label{example1}
R = \left(\begin{array}{cc}
     0 & - x  \\
     - x & - 2y
\end{array}\right).
\end{equation}
The eigenvalues are $\lambda_{1, 2} = - y \pm \sqrt{x^2 + y^2}$. They are not even $C^1$ functions. 

Note, that $R$ is diagonal in every point of the plane, but, at the same time because the eigenfunctions are bad there are no analytic (smooth) diagonalizing coordinate change around the origin. This is the simplest example, that shows, that Haantjes theorem \cite{haantjes} does not work in singular points.
\end{example}


\section{The proof of Theorem \ref{class}}\label{proof1}

We start with lemma
\begin{lemma}\label{lm0}
The algebras in Table 1 and Table 2 in the statement of the Theorem \ref{class} are left-symmetric algebras. Algebras with different names (including different values of the parameters in continuous series) are not isomorphic.
\end{lemma}
{\it Proof. } For algebras $\mathfrak b_{1, \alpha}, \mathfrak b_{2, \beta}, \mathfrak b_3, \mathfrak b_5$, $\mathfrak c_1, \mathfrak c_2, \mathfrak c_3$ and $\mathfrak c_4$ their right-adjoint operators are Nijenhuis by Corollaries \ref{useful} and \ref{useful2}. For $\mathfrak b^+_4, \mathfrak b^-_4, \mathfrak c^+_4, \mathfrak c^-_4$ we check the second condition of Proposition \ref{local1} by direct calculation:
\begin{enumerate}
    \item For $\mathfrak b^+_4$ we have $\operatorname{det} R = x^2, \operatorname{tr} R = - 2y$ and:
\begin{equation*}
    \begin{aligned}
        \left( \begin{array}{cc}
             2x, & 0  \\
        \end{array}\right) = \left( \begin{array}{cc}
             0 ,& - 2  \\
        \end{array}\right)\left( \begin{array}{cc}
            - 2y & x \\
             - x & 0 
        \end{array}\right)
    \end{aligned}
\end{equation*}
    \item For $b^-_4$ we have $\operatorname{det} R = - x^2, \operatorname{tr} R = - 2y$ and:
    \begin{equation*}
    \begin{aligned}
        \left( \begin{array}{cc}
             - 2x, & 0  \\
        \end{array}\right) = \left( \begin{array}{cc}
             0 ,& - 2  \\
        \end{array}\right)\left( \begin{array}{cc}
            - 2y & x \\
             x & 0 
        \end{array}\right)
    \end{aligned}
\end{equation*}
        \item For $\mathfrak c^+_5$ we have $\operatorname{det} R = y^2 - x^2, \operatorname{tr} R = 2y$ and:
        \begin{equation*}
    \begin{aligned}
        \left( \begin{array}{cc}
             - 2x, & 2y  \\
        \end{array}\right) = \left( \begin{array}{cc}
             0 ,& 2  \\
        \end{array}\right)\left( \begin{array}{cc}
            y & - x \\
             - x & y 
        \end{array}\right)
    \end{aligned}
\end{equation*}
        \item For For $\mathfrak c^-_5$ we have $\operatorname{det} R = y^2 + x^2, \operatorname{tr} R = 2y$ and: 
        \begin{equation*}
    \begin{aligned}
        \left( \begin{array}{cc}
             2x, & 2y  \\
        \end{array}\right) = \left( \begin{array}{cc}
             0 ,& 2  \\
        \end{array}\right)\left( \begin{array}{cc}
            y & - x \\
             x & y 
        \end{array}\right)
    \end{aligned}
\end{equation*}
\end{enumerate}
As all right-adjoint operators in Table 1 and Table 2 are Nijenhuis, then by Proposition \ref{main1} all algebras are left-symmetric.

Now, let us proceed with second claim of Lemma. First, note, that algebras from Table 1 and Table 2 are not isomorphic, as their associated Lie algebras are not isomorphic. 

We will need a following observation. Recall, that $\mathfrak a$ has a natural structure of affine manifold. Consider left-symmetric algebras $\mathfrak a, \mathfrak a'$ and their left- and right-adjoint operators $L, R$ and $L', R'$ respectively. We have that $\mathfrak a$ and $\mathfrak a'$ are isomorphic if and only if there exists a linear coordinate change, that transforms $R$ into $R'$ and $L$ into $L'$ simultaneously. Assume now, that there exists a function of four variables $F(y_1, y_2, y_3, y_4)$ such that $F(\operatorname{tr} L, \operatorname{det} L, \operatorname{det} R, \operatorname{tr}R) \equiv 0$
for $\mathfrak a$ and is not identically zero for $\mathfrak a'$. Obviously, coordinate transformations preserve $F$, thus, $\mathfrak a$ and $\mathfrak a'$ are not isomorphic. Now let us proceed with the proof.

Consider $F = \operatorname{det}L$. This function is identically zero on $\mathfrak b_5$ only, thus, this algebra is not isomorphic to $\mathfrak b_{1, \alpha}, \mathfrak b_{2, \beta}, \mathfrak b_3, \mathfrak b_4^+, \mathfrak b_4^-$.

Consider $F = \operatorname{det} R$. It is identically zero for $\mathfrak b_{1, \alpha}$ and $\mathfrak b_3$. Thus, they are not isomorphic to $\mathfrak b_{2, \beta}, \mathfrak b_4^+, \mathfrak b_4^-$. Note, that operator field $L$ for $\mathfrak b_{1, \alpha}$ consists of diagonalizable operators, while $L$ for $\mathfrak b_3$ is Jordan block almost everywhere. Thus, these algebras are not isomorphic.

Fix $\alpha_0$ and consider $F = \alpha_0 \operatorname{det} L - (\operatorname{tr}R)^2$. This function vanishes only for $\mathfrak b_{1, \alpha_0}$ and non-zero for $\alpha \neq \alpha_0$. Thus, algebras $\mathfrak b_{1, \alpha}$ for different values of parameter $\alpha$ are not isomorphic.

Consider $F = (\operatorname{tr} R)^2 - 4 \operatorname{det} R$ (that is the discriminant of characteristic polynomial of $R$). It vanishes identically for $\mathfrak b_{2, \beta}$ and not for $\mathfrak b_4^+, \mathfrak b_4^-$. Thus, $\mathfrak b_{2, \beta}$ is not isomorphic to $\mathfrak b^+_4, \mathfrak b^-_4$.

Fix $\beta_0 \neq 1$ and consider $F = \operatorname{det} L - \big( 1 - \frac{1}{\beta_0}\big) \operatorname{det} R$. It vanishes for $\mathfrak b_{2, \beta_0}$ and does not vanish for $\mathfrak b_{2, \beta}$ for $\beta \neq \beta_0$. Thus, the algebras $\mathfrak b_{2, \beta}$ are not isomorphic for different values of the parameters.

Denote the function
\begin{equation*}
    T(x) = \left\{\begin{array}{cc}
         0 & x \geq 0,  \\
         1 & x < 0,\\
    \end{array}\right.
\end{equation*}
Take $F = T((\operatorname{tr} R)^2 - 4\operatorname{det} R)$. It vanishes for $\mathfrak b_4^-$ and does not vanish for $\mathfrak b_4^+$. This completes the proof for Table 1.

Now consider Table 2. For algebra $\mathfrak c_1$ the operator field $R = L$ is zero, but its not zero for the rest of the Table. Thus, this algebra is not isomorphic to any algebra from the rest of the Table 2.

Take $F = \operatorname{tr} R$. It vanishes for $\mathfrak c_3$, but not for $\mathfrak c_2, \mathfrak c_4, \mathfrak c_5^+, \mathfrak c_5^-$. Thus, $\mathfrak c_3$ is not isomorphic to $\mathfrak c_2, \mathfrak c_4, \mathfrak c_5^+, \mathfrak c_5^-$.

Take $F = \operatorname{det} R$. It vanishes for $\mathfrak c_2$, but not for $\mathfrak c_4, \mathfrak c_5^+, \mathfrak c_5^-$. Thus, the $\mathfrak c_2$ is not isomorphic to $\mathfrak c_4, \mathfrak c_5^+, \mathfrak c_5^-$.

Take $F = T(\operatorname{det} R)$. It is zero for $\mathfrak c_5^-$ and nonzero for $\mathfrak c_5^+$ and $\mathfrak c_4$. Thus, it is not isomorphic to $\mathfrak c_5^+, \mathfrak c_4$. 

Finally, notice, that the right-adjoint operator $R$ for $\mathfrak c_5^+$ almost everywhere has pairwise distinct eigenvalues, while $\mathfrak c_4$ almost everywhere is Jordan block. This completes the proof of Lemma. $\blacksquare$

\begin{lemma}\label{lem1}
Let $R, Q \in \gl(n, \mathbb R)$ and $[R, Q] = \lambda Q$ for $\lambda \neq 0$. Then $Q$ is nilpotent.
\end{lemma}
{\it Proof.} An operator $\mathrm{ad}_R: \gl(n, \mathbb R) \to \gl(n, \mathbb R)$ is defined by the formula $\mathrm{ad}_R Q = [R, Q] = \lambda Q$. From the properties of the matrix commutator it follows that $\mathrm{ad}_R Q^n = n\lambda Q^n$. 

Suppose now that $Q^n \neq 0$ for all $n \in \mathbb N$. It means that the finite dimensional operator $\mathrm{ad}_R$ has an infinite set of eigenvectors $\lambda, 2\lambda, 3\lambda, \dots$. This contradiction completes the proof. $\blacksquare$

\begin{lemma}\label{lem2}
Every two-dimensional commutative subalgebra $\mathfrak h \subset \gl(2, \mathbb R)$ contains the one-dimensional subspace spanned by the identity matrix.
\end{lemma}
{\it Proof}. First, recall some facts about $\gl(n, \mathbb R)$. It has an invariant scalar product, which identifies $\gl(n, \mathbb R)$ with its dual. In this construction the adjoint orbits are identified with coajoint orbits. The latter are simplectic manifolds and, thus, even dimensional. 

At the same time the dimension of the adjoint orbit $\mathcal O_X$ of a given element $X$ is
$$
\operatorname{dim} \mathcal O_X =\operatorname{dim} \gl(n, \mathbb R) - \operatorname{dim} \mathfrak z_X,
$$
where $\mathfrak z_X$ denotes the centralizer of $X$, that is all elements $Y \in \gl(n, \mathbb R)$, that commute with $X$. For $n = 2$ we get that $4 - \operatorname{dim} \mathfrak z_X$ must be even, thus, the dimension of centralizer of any element $X$ is even.

Now consider element $X \in \mathfrak h$. By definition $\mathfrak h \subseteq \mathfrak z_X$. This implies, that dimension of $\mathfrak z_X$ is either $2$ or $4$. If $\operatorname{dim} \mathfrak z_X = 4$, then $X$ lies in the center of $\gl(2, \mathbb R)$ and is proportional to identity matrix. Thus, the statement of the Lemma holds in thic case. If $\operatorname{dim} \mathfrak z_X = 2$, then $\mathfrak h$ coincides with $\mathfrak z_X$. In particular identity matrix commutes with $X$ and lies in $\mathfrak h$. The Lemma is proved. $\blacksquare$.

For the rest of the proof let us forget about affine structure on $\mathfrak a$ and treat $L$ simply as left-adjoint action of $\mathfrak a$ on itself. Fixing basis $\eta_1, \eta_2$ and we get, that $L_{\xi}$ defines a map from $\mathfrak a$ to matrix algebra $\mathfrak{gl}(2, \mathbb R)$. Identity \eqref{left} implies, that the image of this map is a subalgebra in $\gl(2, \mathbb R)$ (in fact this is the image of the adjoint action of associated Lie algebra). We denote the dimension of this subalgebra as $N$.

\subsection{$N = 0$}
This implies that $L_{\eta} = 0$ for an arbitrary $\eta \in \mathfrak a$. Thus, $\mathfrak a$ is isomorphic to $\mathfrak c_1$ from Table 2.

\subsection{$N = 1$}

W.l.o.g we may assume, that in basis $\xi_1, \xi_2$ the element $\xi_1$ spans the kernel of map $L$, that is $\xi_1 \star \xi_1 = 0$, $\xi_1 \star \xi_2 = 0$ or, in other words, $L_{\xi_1} = 0$. From identity \eqref{left} we get $L_{[\xi_1, \xi_2]} = [L_{\xi_1}, L_{\xi_2}] = 0$. This implies, that
\begin{equation}\label{gamma}
[\xi_1, \xi_2] = \gamma \xi_1.    
\end{equation}

If $\gamma = 0$ in \eqref{gamma} we get, that $[\xi_1, \xi_2] = 0$. Thus, in given basis structure relations for $\mathfrak a$ are
\begin{equation}\label{pref0}
\begin{aligned}
& \xi_1\star\xi_1 = \xi_1\star\xi_2 = \xi_2\star\xi_1 = 0, \\
& \xi_2\star\xi_2 = b \xi_1 + a \xi_2 , \\    
\end{aligned}
\end{equation}
where $a, b$ some constants.

If $a \neq 0$ in \eqref{pref0}, consider the change of basis $\eta_1 = \xi_1$ and $\eta_2 = \frac{b}{a^2} \xi_1 + \frac{1}{a} \xi_2$. The structure relations \eqref{pref0} take form
\begin{equation*}
\begin{aligned}
& \eta_1\star\eta_1 = \eta_1\star\eta_2 = \eta_2\star\eta_1 = 0, \\
& \eta_2\star\eta_2 = \frac{1}{a^2} \xi_2 \star \xi_2 = \Big(\frac{b}{a^2} \xi_1 + \frac{1}{a} \xi_2\Big) = \eta_2. \\    
\end{aligned}
\end{equation*}
This yields the structure relations of the left-symmetric algebra $\mathfrak c_2$ from Table 2.

If $a = 0$ in \eqref{pref0}, then $b \neq 0$ (otherwise $\operatorname{dim}\operatorname{Im} L = 0$). After the change of basis $\eta_1 = b \xi_1, \eta_2 = \xi_2$ relations \eqref{pref0} take form
\begin{equation*}
\begin{aligned}
& \eta_1\star\eta_1 = \eta_1\star\eta_2 = \eta_2\star\eta_1 = 0, \\
& \eta_2\star\eta_2 = b \xi_1 = \eta_1. \\    
\end{aligned}
\end{equation*}
This yields the structure relations of the left-symmetric algebra $\mathfrak c_3$ from Table 2.

Assume now, that $\gamma \neq 0$ in \eqref{gamma}. We have $[\xi_1, \xi_2] = \gamma \xi_1$. Changing coordinates as $\xi'_1 = \xi_1$ and $\xi'_2 = - \frac{1}{\gamma} \xi_2$ we get, that $[\xi'_1, \xi'_2] = - \xi'_1$. So, w.l.o.g. we may assume that if $\gamma \neq 0$ in \eqref{gamma}, then $\gamma = - 1$. 

The latter implies, that $\xi_1 \star \xi_2 - \xi_2 \star \xi_1 = \xi_1$. The structure relations in this case are 
\begin{equation}\label{pref}
\begin{aligned}
    & \xi_1\star\xi_1 = \xi_1\star\xi_2 = 0, \\
    & \xi_2\star\xi_1 = \xi_1, \\
    & \xi_2\star\xi_2 = b \xi_1 + a \xi_2,
\end{aligned}    
\end{equation}
for some constants $a$ and $b$.

If $a \neq 1$ in \eqref{pref}, then after the change of basis $\eta_1 = \xi_1, \eta_2 = - \frac{b}{1 - a} \xi_1 + \xi_2$ the relations \eqref{pref} take form
\begin{equation*}
\begin{aligned}
    & \eta_1 \star \eta_1 = \eta_1 \star \eta_2 = 0, \\
    & \eta_2 \star \eta_1 = \Big( - \frac{b}{1 - a} \xi_1 + \xi_2\Big) \star \xi_1 = \xi_1 = \eta_1, \\
    & \eta_2 \star \eta_2 = \Big( - \frac{b}{1 - a} \xi_1 + \xi_2\Big) \star \Big( - \frac{b}{1 - a} \xi_1 + \xi_2\Big) = - \frac{b}{1 - a} \xi_1 + b\xi_1 + a\xi_2 = \\
    & = a \Big( - \frac{b}{1 - a} \xi_1 + \xi_2\Big) = a \eta_2.
\end{aligned}
\end{equation*}
Renaming $a$ as $\alpha$, we get the structure relations of the left-symmetric algebra $\mathfrak b_{1, \alpha}$ from Table 1.

If in \eqref{pref} constants $a = 1$ and $b = 0$, then we get
\begin{equation*}
    \begin{aligned}
    & \xi_1\star\xi_1 = \xi_1\star\xi_2 = 0, \\
    & \xi_2\star\xi_1 = \xi_1, \\
    & \xi_2\star\xi_2 = \xi_2,
\end{aligned}   
\end{equation*}
These are the structure relations of $\mathfrak b_{1, 1}$ from Table 1.

Assume, finally, that in \eqref{pref} $a = 1$ and $b \neq 0$. After the change of basis $\eta_1 = b \xi_1, \eta_2 = \xi_2$ the relations \eqref{pref} take the form
\begin{equation*}
    \begin{aligned}
        & \eta_1 \star \eta_1 = \eta_1 \star \eta_2 = 0, \\
        & \eta_2 \star \eta_1 = b \xi_2 \star \xi_1 = b \xi_1 = \eta_1, \\
        & \eta_2 \star \eta_2  = \xi_2 \star \xi_2 = b \xi_1 + \xi_2 = \eta_1 + \eta_2. \\
    \end{aligned}
\end{equation*}
These are the structure relations of $\mathfrak b_5$ from Table 1.

\subsection{$N = 2$ and associated Lie algebra is abelian}

$L$ defines a faithful representation of the associated Lie algebra for $\mathfrak a$. As associated Lie algebra is abelian, then by Lemma \ref{lem2}, the image of the representation contains the identity matrix $\operatorname{Id}$. W.l.o.g. we may assume that $L_{\xi_2} = \operatorname{Id}$ in given basis $\xi_1, \xi_2$. In other words $\xi_1 \star \xi_1 = \xi_1$ and $\xi_1 \star \xi_2 = \xi_2$.

Identity \eqref{left} implies, that $L_{[\xi_1, \xi_2]} = [L_{\xi_1}, L_{\xi_2}] = [L_{\xi_1}, \operatorname{Id}] = 0$. As representation is faithful and associated Lie algebra is abelian we have $[\xi_1, \xi_2] = 0$. The structure relations in this case are:
\begin{equation}\label{pref1}
    \begin{aligned}
        & \xi_2 \star \xi_1 = \xi_1 \star \xi_2 = \xi_1, \\
        & \xi_2 \star \xi_2 = \xi_2, \\
        & \xi_1\star\xi_1 = a \xi_1 + b \xi_2, \\
    \end{aligned}
\end{equation}
for some constants $a$ and $b$.

Assume that in \eqref{pref1} coefficients satisfy the relation $\frac{a^2}{4} + b = 0$. Then after the change of basis $\eta_1 = \xi_1 - \frac{a}{2}\xi_2, \eta_2 = \xi_2$ the relations \eqref{pref1} take form
\begin{equation*}
    \begin{aligned}
        & \eta_2 \star \eta_1 = \eta_1, \quad \eta_1 \star \eta_2 = \Big(\xi_1 - \frac{a}{2}\xi_2\Big) \star \xi_2 = \xi_1 - \frac{a}{2}\xi_2 = \eta_1, \\
        & \eta_2 \star \eta_2 = \eta_2, \quad \Big(\xi_1 - \frac{a}{2}\xi_2\Big) \star \Big( \xi_1 - \frac{a}{2}\xi_2 \Big) = a \xi_1 + b \xi_2 - a \xi_1 + \frac{a^2}{4} \xi_2 = 0.
    \end{aligned}
\end{equation*}
These are the structure relations of $\mathfrak c_4$ from Table 2.

Assume, that the coefficients $a, b$ in \eqref{pref1} are such, that $\frac{a^2}{4} + b \neq 0$. After the change of basis $\eta_1 = \frac{1}{\sqrt{\vert \frac{a^2}{4} + b}\vert}\xi_1 - \frac{a}{2\sqrt{\vert\frac{a^2}{4} + b}\vert}\xi_2, \eta_2 = \xi_2$ the relations \eqref{pref1} take the form
\begin{equation*}
    \begin{aligned}
        & \eta_2 \star \eta_1 = \eta_1, \quad \eta_1 \star \eta_2 = \eta_1, \quad \eta_2 \star \eta_2 = \eta_2, \\
        & \eta_1 \star \eta_1 = \Big(\frac{1}{\sqrt{\vert \frac{a^2}{4} + b}\vert}\xi_1 - \frac{a}{2\sqrt{\vert\frac{a^2}{4} + b}\vert}\xi_2\Big) \star \Big( \frac{1}{\sqrt{\vert \frac{a^2}{4} + b}\vert}\xi_1 - \frac{a}{2\sqrt{\vert\frac{a^2}{4} + b}\vert}\xi_2 \Big) = \frac{\frac{a^2}{4} + b}{\vert \frac{a^2}{4} + b \vert} \eta_2.
    \end{aligned}
\end{equation*}
Depending on sign of $\frac{a^2}{4} + b$ we get the structure relations of either $c^+_5$ or $c^-_5$.

\subsection{$N = 2$ and associated Lie algebra is non-abelian}

As associated Lie algebra of $\mathfrak a$ is not abelian, then w.l.o.g. we may assume that $[\xi_1, \xi_2] = \xi_1\star\xi_2 - \xi_2\star\xi_1 = \xi_1$. Identity \eqref{left} implies, that $[L_{\xi_1}, L_{\xi_2}] = L_{\xi_1}$. By Lemma \ref{lem1} $L_{\xi_1}$ is nilpotent and, as $L_{\xi_1} \neq 0$ as $L_{\xi}$ is a faithful representation of associated Lie algebra. In case of dimension two, the kernel and the image of $L_{\xi_1}$ coincide. We have two cases.

\subsubsection{Image of $L_{\xi_1}$ is spanned by $\xi_1$}

The condition implies, that $L_{\xi_1} \xi_1 = \xi_1 \star \xi_1 = 0$ and $L_{\xi_1} \xi_2 = \xi_1 \star \xi_2 = \gamma \xi_1$ with $\gamma \neq 0$. Assume $\xi_2 \star \xi_2 = a \xi_1 + b \xi_2$. The left-symmetry of associator yields:
\begin{equation*}
    \begin{aligned}
       0 = & \mathcal A(\xi_1, \xi_2, \xi_2) - \mathcal A(\xi_2, \xi_1, \xi_2) = \\
       = &(\xi_1 \star \xi_2) \star \xi_2 - \xi_1 \star (\xi_2 \star \xi_2) - (\xi_2 \star \xi_1) \star \xi_2 + \xi_2 \star (\xi_1 \star \xi_2) = \\
        = & \gamma^2 \xi_1 - b\gamma \xi_1 - \gamma (\gamma - 1)\xi_1 + \gamma (\gamma - 1) \xi_1 = \\
        = & \gamma (b - \gamma)\xi_1.
    \end{aligned}
\end{equation*}
Thus, we get the following structure relations
\begin{equation}\label{pref2}
    \begin{aligned}
        & \xi_1 \star \xi_1 = 0, \quad \xi_1 \star \xi_2 = \gamma \xi_1, \\
        & \xi_2 \star \xi_1 = (\gamma - 1)\xi_1, \quad \xi_2 \star \xi_2 = a \xi_1 + \gamma \xi_2.
    \end{aligned}
\end{equation}

Assume, that in \eqref{pref2} $\gamma = 1$ and $a \neq 0$. After the change of basis $\eta_1 = a \xi_1, \eta_2 = \xi_2$ the structure relations \eqref{pref2} take form
\begin{equation*}
    \begin{aligned}
        & \eta_1 \star \eta_1 = 0, \quad \eta_1 \star \eta_2 = \eta_1, \\
        & \eta_2 \star \eta_1 = 0, \quad \eta_2\star \eta_2 = \eta_1 + \eta_2.
    \end{aligned}
\end{equation*}
These are the structure relations of $\mathfrak b_5$.

Assume, that in \eqref{pref2} $\gamma = 1$ and $a = 0$. Renaming the basis $\eta_1 = \xi_1, \eta_2 = \xi_2$ we get the structure relations
\begin{equation*}
    \begin{aligned}
        & \eta_1 \star \eta_1 = 0, \quad \eta_1 \star \eta_2 = \eta_1, \\
        & \eta_2 \star \eta_1 = 0, \quad \eta_2 \star \eta_2 = \eta_2.
    \end{aligned}
\end{equation*}
These are the structure relations of $\mathfrak b_{2, 1}$.

Assume, that in \eqref{pref2} $\gamma \neq 1$. After change of basis $\eta_1 = \xi_1, \eta_2 = \frac{a}{\gamma (1 - \gamma)}\xi_1 + \frac{1}{\gamma} \xi_2$ the relations \eqref{pref2} take the form 
\begin{equation*}
    \begin{aligned}
        & \eta_1 \star \eta_1 = 0, \quad \eta_1 \star \eta_2 = \xi_1 \star \Big( \frac{a}{\gamma (1 - \gamma)}\xi_1 + \frac{1}{\gamma} \xi_2\Big) = \xi_1 = \eta_1, \\
        & \eta_2 \star \eta_1 = \Big( \frac{a}{\gamma (1 - \gamma)}\xi_1 + \frac{1}{\gamma} \xi_2\Big) \star \eta_1 = \frac{\gamma - 1}{\gamma} \xi_1 = \Big( 1 - \frac{1}{\gamma}\Big) \eta_1, \\
        & \eta_2 \star \eta_2 = \Big( \frac{a}{\gamma (1 - \gamma)}\xi_1 + \frac{1}{\gamma} \xi_2\Big) \star \Big( \frac{a}{\gamma (1 - \gamma)}\xi_1 + \frac{1}{\gamma} \xi_2\Big) = \\
        & = - \frac{a}{\gamma^2} \xi_1 + \frac{a}{\gamma (1 - \gamma)} \xi_1 + \frac{a}{\gamma^2} \xi_1 + \frac{1}{ \gamma} \xi_2 = \eta_2.
    \end{aligned}
\end{equation*}
Renaming $\gamma$ as $\beta$ we obtain the structure relations for $\mathfrak b_{2, \beta}$.

\subsubsection{Image of $L_{\xi_1}$ is spanned by $a \xi_1 + \xi_2$ for some constant $a$} 

Recall, that in dimension two the image and kernel of nilpotent operator coincide. As $\xi_1$ and $a \xi_1 + \xi_2$ are linearly independent for all $a \in \mathbb R$, we have $L_{\xi_1} \xi_1 = \xi_1 \star \xi_1 = b (a \xi_1 + \xi_2)$ for some $b \neq 0$. Consider change of basis $\eta_1 = \frac{1}{\sqrt{|b|}}\xi_1, \eta_2 = \xi_2 + a \xi_1$. By definition we have
$$
L_{\eta_1} \eta_1 = \eta_1 \star \eta_1 = \frac{b}{|b|} (a\xi_1 + \xi_2) = \operatorname{sgn} (b) \, \eta_2.
$$
At the same time by definition $L_{\eta_1} \eta_2 = \eta_1 \star \eta_2 = 0$ and
\begin{equation*}
    \begin{aligned}{}
        [\eta_1, \eta_2] = \frac{1}{\sqrt{|b|}}[\xi_1, \xi_2] = \frac{1}{\sqrt{|b|}} \xi_1 = \eta_1.
    \end{aligned}
\end{equation*}
The left symmetry of the associator for triple $\eta_1, \eta_2, \eta_2$ yields:
\begin{equation*}
    \begin{aligned}
       0 = & \mathcal A(\eta_1, \eta_2, \eta_2) - \mathcal A(\eta_2, \eta_1, \eta_2) = \\
       = &(\eta_1 \star \eta_2) \star \eta_2 - \eta_1 \star (\eta_2 \star \eta_2) - (\eta_2 \star \eta_1) \star \eta_2 + \eta_2 \star (\eta_1 \star \eta_2) = \\
        = & - \eta_1 \star (\eta_2 \star \eta_2). 
    \end{aligned}
\end{equation*}
This implies, that $\eta_2 \star \eta_2 = \gamma \eta_2$ for some constant $\gamma$. The left symmetry of the associator for triple $\eta_1, \eta_2, \eta_1$ yields:
\begin{equation*}
    \begin{aligned}
        0 = & \mathcal A(\eta_1, \eta_2, \eta_1) - \mathcal A(\eta_2, \eta_1, \eta_1) = \\
        = & (\eta_1 \star \eta_2) \star \eta_1 - \eta_1 \star (\eta_2 \star \eta_1) - (\eta_2 \star \eta_1)\star \eta_1 + \eta_2 \star (\eta_1 \star \eta_1) = \\
        = & 2 \eta_1 \star \eta_1 + \operatorname{sgn}(b) \eta_2 \star \eta_2 =  \operatorname{sgn}(b) (2 + \gamma) \eta_2.
    \end{aligned}
\end{equation*}
As $\operatorname{sgn}(b) \neq 0$, then $\gamma = - 2$. Thus, we get the following structure relations:
\begin{equation*}
    \begin{aligned}
        & \eta_1 \star \eta_1 = \operatorname{sgn} (b) \, \eta_2, \quad \eta_1 \star \eta_2 = 0, \\
        & \eta_2 \star \eta_1 = - \eta_1, \quad \eta_2 \star \eta_2 = - 2 \eta_2.
    \end{aligned}
\end{equation*}
Depending on the sign of $b$ this yields the structure relations for either $\mathfrak b^+_4$ or $\mathfrak b^-_4$. The Theorem is proved. $\blacksquare$




\section{Proof of Proposition \ref{main2}}\label{proof2}

Let us recall that we have a pair of vector fields $v, w$ with the property $v({\mathrm{p}}) = v_{\mathrm{p}}$ and $w({\mathrm{p}}) = w_{\mathrm{p}}$, and the operation 
$$
v_{\mathrm{p}} \star w_{\mathrm{p}} = ([Rv, w] - R[v, w]) \vert_{\mathrm{p}}. 
$$
In local coordinates $x^1, \dots, x^n$ we have (the underlined terms cancel):
\begin{equation*}
    \begin{aligned}
        & \big( [Rv, w] - R[v, w]\big)^k = \\
        = & \pd{(R^k_{\alpha} v^{\alpha})}{x^{\beta}} w^{\beta} - \pd{w^k}{x^{\beta}} R^{\beta}_{\alpha} v^{\alpha} - R^k_{\alpha} \pd{v^{\alpha}}{x^{\beta}} w^{\beta} + R^k_{\alpha} \pd{w^{\alpha}}{x^{\beta}} v^{\beta} = \\
        = & \pd{R^k_{\alpha}}{x^{\beta}} v^{\alpha} w^{\beta} + \underline{R^k_{\alpha} \pd{v^{\alpha}}{x^{\beta}} w^{\beta}} - \pd{w^k}{x^{\beta}} R^{\beta}_{\alpha} v^{\alpha} - \underline{R^k_{\alpha} \pd{v^{\alpha}}{x^{\beta}} w^{\beta}} + R^k_{\alpha} \pd{w^{\alpha}}{x^{\beta}} v^{\beta} = \\
        = & \pd{R^k_{\alpha}}{x^{\beta}} v^{\alpha} w^{\beta} + R^k_{\alpha} \pd{w^{\alpha}}{x^{\beta}} v^{\beta} - \pd{w^k}{x^{\beta}} R^{\beta}_{\alpha} v^{\alpha}.
    \end{aligned}
\end{equation*}
Substituting point ${\mathrm{p}}$, where $R^k_i = \lambda \delta^k_i$ into the last equation, we get
$$
\big( [Rv, w] - R[v, w]\big)^k \vert_{\mathrm{p}} = \pd{R^k_{\alpha}}{x^{\beta}} \vert_{\mathrm{p}} \, v_{\mathrm{p}}^{\alpha} w_{\mathrm{p}}^{\beta}.
$$
Thus, we see, that the definition of the operation for $v_{\mathrm{p}}, w_{\mathrm{p}}$ does not depend on the continuations $v, w$ and that the structure constants of the algebra are $\pd{R^k_{\alpha}}{x^{\beta}} \vert_{\mathrm{p}}$. 

Let us now show, that the algebra is left-symmetric. Denote $a^k_{ij} = \pd{R^k_i}{x^j}\vert_{\mathrm{p}}$. Recall, that if $R$ is Nijenhuis, then $R - \lambda \operatorname{Id}$ is Nijenhuis. So w.l.o.g. assume that $R = 0$ at point ${\mathrm{p}}$ and consider
\begin{equation*}
    \begin{aligned}
        0 = & \pd{}{x^r} \Big( \mathcal N_R \Big) \vert_{\mathrm{p}} = \\
        = & \pd{}{x^r} \Bigg(\pd{R^{\alpha}_i}{x^j} R_{\alpha}^k - \pd{R^{\alpha}_j}{x^i} R_{\alpha}^k -\pd{R^k_i}{x^{\alpha}} R^{\alpha}_j + \pd{R^k_j}{x^{\alpha}} R^{\alpha}_i\Bigg)\vert_{\mathrm{p}} = \\
        = & a^{\alpha}_{ij} a_{\alpha r}^k - a^{\alpha}_{ji} a_{\alpha r}^k - a^k_{i \alpha} a^{\alpha}_{jr} + a^k_{j \alpha} a^{\alpha}_{ir}.
    \end{aligned}
\end{equation*}
In natural basis $\eta_1, \dots, \eta_n$ in $T_{\mathrm{p}} M^n$, associated with coordinates $x^1, \dots, x^n$ on the neighbourhood of $\mathrm{p}$, the last equation yields
\begin{equation*}
\begin{aligned}
0 = (\eta_i \star \eta_j) \star \eta_r - (\eta_j \star \eta_i)\star \eta_r - & \eta_i \star (\eta_j \star \eta_r) + \eta_j \star (\eta_i \star \eta_r) = \\
= & \mathcal A(\eta_i, \eta_j, \eta_r) - \mathcal A(\eta_j, \eta_i, \eta_r).
\end{aligned}
\end{equation*}
Thus, vanishing of the Nijenhuis torsion implies, that algebra $\mathfrak a$ is left-symmetric and Proposition \ref{main2} is proved. $\blacksquare$


\section{Proof of Theorem \ref{themain1} and Theorem \ref{themain2}}\label{proof3}

\subsection{Algebras $\mathfrak c_1, \mathfrak c_2, \mathfrak c_3, \mathfrak c_4$}\label{romm}

For each of these four left-symmetric algebras we provide polynomial Nijenhuis operators $R$ with $R_1$ coinciding with the right-adjoint operator of the corresponding algebra from Table 2 of Theorem \ref{class}. We also give polynomial function $F (\operatorname{tr} R, \operatorname{det} R)$ that identically vanish for $R_1$ but non-zero on every neighbourhood of the coordinate origin. This implies, that there are no linearizing coordinate change for $R$ in both analytic and smooth category.

\begin{equation*}
    \begin{aligned}
      \mathfrak c_1: & \quad  R = \left(\begin{array}{cc}
x^2 & 0 \\
0 & y^2 \\
\end{array}\right), \quad F = \operatorname{tr} R = x^2 + y^2; \\
\mathfrak c_2: &  \quad R = \left(\begin{array}{cc}
x^2 & 0 \\
0 & y \\
\end{array}\right), \quad F = \operatorname{det} R = x^2y; \\
\mathfrak c_3: & \quad R = \left(\begin{array}{cc}
y^2 & y \\
0 & y^2 \\
\end{array}\right), \quad F = \operatorname{det} R = y^4;\\
\mathfrak c_4: & \quad R = \left(\begin{array}{cc}
y + yx^2 & x + x^3\\
-xy^2 & y - yx^2\\
\end{array}\right), \quad F = (\operatorname{tr} R)^2 - 4 \operatorname{det} R = - 4 x^2y^2.
    \end{aligned}
\end{equation*}

The operators $R$ for $\mathfrak c_1, \mathfrak c_2$ are Nijenhuis by Haantjes theorem \cite{haantjes}, $R$ for $\mathfrak c_3$ is Nijenhuis by Corollary \ref{useful}. Finally operator $R$ for $\mathfrak c_4$ is Nijenhuis by formula \ref{local2} for $\operatorname{tr}R = 2y$ and $\operatorname{det} R = y^2 + x^2 y^2$.

\subsection{Algebra $\mathfrak b_5$}

Consider operator field
$$
R = \left(\begin{array}{cc}
y & y - x^2 \\
0 & y \\
\end{array}\right).
$$
Its linear part coincides with the right-adjoint operator of $\mathfrak b_5$ from Table 1 of Theorem \ref{class}. At the same time the curve $y = x^2$ through the coordinate origin for $y \neq 0$ consists of operators, proportional to $\operatorname{Id}$ with non-zero coefficient. Thus, in every neighbourhood of the coordinate origin there exist a point, for which $R$ is proportional to the identity with some non-zero coefficient.

At the same time the algebraic type of the right-adjoint operator of $\mathfrak b_5$ in every point is either Jordan block or zero operator. Thus, there are no analytic or smooth linearizing coordinate change.

\subsection{Algebras $\mathfrak c_5^+, \mathfrak c_5^-$}

Consider two dimensional real plane with coordinates $x, y$ and Nijenhuis operator field $R$, vanishing at the origin. We assume, that the linear part of Taylor expansion of $R$ at the origin is the right-adjoint operator of $\mathfrak c^+_5$ from Table 2 of Theorem \ref{class}. 

This implies, that the Taylor expansion of $\operatorname{tr} R$ at the origin starts with linear term $2y$ and the Taylor expansion of $g = \operatorname{det} R$ starts with quadratic term $y^2 - x^2$. In particular, $\frac{\partial^2 g}{\partial x^2} (0, 0) = -2 \neq 0$, $\frac{\partial^2 g}{\partial x \partial y}(0, 0) = 0$ and $\frac{\partial^2 g}{\partial y^2} (0, 0) = 2$. W.l.o.g. assume, that $\operatorname{tr} R = 2y$. 

Applying Lemma \ref{morse} we get that there exist new coordinates $\bar x, \bar y$, such that
$$
\operatorname{tr} R = 2 \bar y, \quad g = - \bar x^2 + k(\bar y)
$$
and $k''(0) = \frac{\partial^2 g}{\partial y^2} (0, 0) = 2$. Substituting this expression into the condition \eqref{condition} and taking $x = 0$ we get
$$
\frac{1}{2}k'(\bar y) (2\bar y - \frac{1}{2}k'(\bar y)) - k(\bar y) = 0.
$$
Differentiating both sides by $\bar y$ we get
$$
\frac{1}{2}k''(\bar y) (2 \bar y - k') = 0.
$$
As $k(0) = 0$ and $k''(0) = 2$, we have that $k' = 2 \bar y$ and, thus, $g = \operatorname{det} R =  \bar y^2 - \bar x^2$. Applying formula \eqref{local2} to the given trace and determinant, we get, that in coordinates $\bar x, \bar y$ Nijenhuis operator $R$ coincides with the right-adjoint operator of $\mathfrak c^+_5$. The proof for $\mathfrak c^-_5$ is similar and both left-symmetric algebras are non-degenerate in smooth and analytic category.

\subsection{Algebras $\mathfrak b^+_4, \mathfrak b^-_4$}

The proof for $\mathfrak b^+_4$ follows the same steps as the proof for $\mathfrak c^+_5$ until the equation
$$
\frac{1}{2}k''(\bar y) (2 \bar y - k') = 0
$$
is obtained. We have, that $k''(0) = \frac{\partial^2 g}{\partial y^2} = 0$. Denote $F = 2 \bar y - k'$ and $F'(0) = 2$, thus, function $F \neq 0$ for sufficiently small neighbourhood of $y = 0$. This implies that $k''$ identically vanish and we get $g = \operatorname{der} R = \bar x^2$. Applying formula \eqref{local2} we get, that Nijenhuis operator $R$ in coordinates $\bar x, \bar y$ is the right-adjoint operator of $\mathfrak b^+_5$. The proof for $\mathfrak b^+_4$ is similar. Thus, both $\mathfrak b^+_4$ and $\mathfrak b^-_4$ are non-degenerate in smooth and analytic category.

\subsection{Algebra $\mathfrak b_{2, \beta}$}

Consider two dimensional real plane with coordinates $x, y$ and Nijenhuis operator $R$ in the form
\begin{equation}\label{fm}
    R = \left(\begin{array}{cc}
         y & f(x, y)  \\
         0 & y 
    \end{array}\right).
\end{equation}
Assume that $R$ vanishes at the coordinate origin and linear part of its Taylor expansion is the right-adjoint operator of $\mathfrak b_{2, \beta}$ from Table 1 of Theorem \ref{class}. In particular, the linear part of Taylor expansion of $f(x, y)$ is $\big(1 - \frac{1}{\beta}\big)x$. We have two cases.

\subsubsection{Case $\beta \neq 1$}

In this case $\pd{f}{x} (0, 0) \neq 0$. By the Implicit Function Theorem there exists a curve $r(y)$, such that $f(r(y), y) = 0$ and $r(0) = 0$. This is exactly the curve, that consists of operators $R$, that are multiple of $\operatorname{Id}$. We call this curve \textbf{the characteristic curve}. Next lemma provides necessary condition for the linearization in case $\beta \neq 1$ for both analytic and smooth category.

\begin{lemma}
Consider Nijenhuis analytic (smooth) operator $R$ in the form \eqref{fm}. Assume, that there exists an analytic (smooth) linearizing coordinate change and denote the characteristic curve as $r(y)$. Then $\pd{f}{x} (r(y), y)$ is constant function on sufficiently small neighbourhood of $y = 0$.
\end{lemma}\label{necessary}
{\it Proof. }Consider linearizing coordinate change $\bar x = g(x, y), \bar y = h(x, y)$. By definition $\bar y = y = \frac{1}{2} \operatorname{tr} R$. Thus, the linearing coordinate change is in triangular form $\bar x = g(x, y), \bar y = y$. Operator field $R$ transformed as
\begin{equation*}
    R = \left(\begin{array}{cc}
         \pd{g}{x} & \pd{g}{y}  \\
          0 & 1 
    \end{array}\right) \left( \begin{array}{cc}
         y & \Big(1 - \frac{1}{\beta}\Big)g  \\
         & y
    \end{array}\right) \left(\begin{array}{cc}
         \frac{1}{\pd{g}{x}} & - \frac{\pd{g}{y}}{\pd{g}{x}} \\
          0 & 1 
    \end{array}\right) = \left( \begin{array}{cc}
         y &  \Big(1 - \frac{1}{\beta}\Big) \frac{g}{\pd{g}{x}}\\
         0 & y  
    \end{array}\right)
\end{equation*}
We have
$$
f(x, y) = \big(1 - \frac{1}{\beta}\big)\frac{g}{\pd{g}{x}}.
$$
Note, that as $\pd{g}{x} (0, 0) \neq 0$, thus, for $r(y)$ we have $g(r(y), y) = 0$ for sufficiently small neighbourhood of $y = 0$. We get
\begin{equation*}
    \pd{f}{x} = \big(1 - \frac{1}{\beta}\big)\Big( 1 - g \frac{1}{\big(\pd{g}{x}\big)^2} \frac{\partial^2 g}{\partial x^2}\Big).
\end{equation*}
Substituting $x = r(y)$ we get, that r.h.s. of this equation identically equals to $(1 - \frac{1}{\beta})$. $\blacksquare$

Consider polynomial operator field 
\begin{equation*}
    R = \left(\begin{array}{cc}
         y & \Big(1 - \frac{1}{\beta}\Big)x + x^2 + xy  \\
         0 & y 
    \end{array}
    \right)
\end{equation*}
By Corollary \ref{useful} $R$ is Nijenhuis and it is in the form \eqref{fm}. For sufficiently small $x, y$ the equation $\Big(\Big(1 - \frac{1}{\beta}\Big) + x + y \Big)x = 0$ defines the curve $x = 0$. Thus, this is the characteristic curve $r(y) \equiv 0$. At the same time $\pd{f}{x} (r(y), y) = 1 - \frac{1}{\beta} + y$, which is non-constant. Thus, by Lemma \ref{necessary} $R$ does not satisfy the necessary condition for linearization and has no linearizing coordinate change for both analytic and smooth case.

\subsubsection{Case $\beta = 1$}

In this case all right-adjoint operators of $\mathfrak b_{2, 1}$ are the multiples of the $\operatorname{Id}$. Consider polynomial operator field 
\begin{equation*}
    R = \left(\begin{array}{cc}
         y & x^2 + y^2  \\
         0 & y 
    \end{array}
    \right)
\end{equation*}
By Corollary \ref{useful} $R$ is Nijenhuis and the linear part in its Taylor expansion coincides with the right-adjoint operator of $\mathfrak b_{2, 1}$. At the same time for $x, y \neq 0$ the algebraic type $R$ is Jordan block, while the right-adjoint operator of $\mathfrak b_{2, 1}$ is diagonalizable at every point. Thus, $\mathfrak b_{1, 1}$ is degenerate in both analytic and smooth category.

\subsection{Algebra $\mathfrak b_{1, \alpha}$ (analytic case)}

We start with lemma.

\begin{lemma}\label{prenormal}
Consider analytic (smooth) Nijenhuis operator $R$ on the real plane with coordinates $x, y$ that vanishes at the origin. Assume, that $R$ is in the form
\begin{equation*}
    R = \left(\begin{array}{cc}
         0 & f(x, y)  \\
         0 & \alpha y 
    \end{array}
    \right),
\end{equation*} for some constant $\alpha$ and that the linear part of the Taylor expansion of $f(x, y)$ is $x + \beta y$. There exists a analytic (smooth) linearizing coordinate change, that transforms $R$ into
\begin{equation*}
    R = \left(\begin{array}{cc}
         0 & \bar x + \beta \bar y   \\
         0 & \alpha \bar y 
    \end{array}
    \right)
\end{equation*}
if and only if there exists a linearizing coordinate change for vector field $v = (f(x, y), y)^T$.
\end{lemma}
{\it Proof. } First, assume that there exists a linearizing coordinate change $\bar x = g(x, y), \bar y = h(x, y)$ for Nijenhuis operator $R$. Note, that $\operatorname{tr} R = \alpha y = \alpha \bar y$, thus, the coordinate change has the " triangular" form $\bar x = g(x, y), \bar y = y$. Operator field transformed as
\begin{equation}\label{swipe}
\begin{aligned}
     \left(\begin{array}{cc}
        \pd{g}{x} & \pd{g}{y}  \\
         0 & 1 
    \end{array}\right)
    \left(\begin{array}{cc}
0 & f(x, y) \\
0 & \alpha y \\
\end{array}\right)& \left(\begin{array}{cc}
        \frac{1}{\pd{g}{x}} & - \frac{\pd{g}{y}}{\pd{g}{x}}  \\
         0 & 1 
    \end{array}\right) = \left(\begin{array}{cc}
0 & \pd{g}{x} f(x, y) + \alpha \pd{g}{y} y  \\
0 & \alpha y \\
\end{array}\right) = \\
& = \left(\begin{array}{cc}
        0 & g (x, y) + \beta y  \\
         0 & \alpha y 
    \end{array}\right) = \left(\begin{array}{cc}
        0 & \bar x + \beta \bar y  \\
         0 & \alpha \bar y 
    \end{array}\right),
\end{aligned}
\end{equation}
As $\pd{g}{x} f(x, y) + \alpha \pd{g}{y} y = v(g)$ this can be written as $v(y) = \alpha y, v(g) = g + \beta y$. Thus, the analytic (smooth) linearizing coordinate change $\bar x = x, \bar y = y$ for $R$ is at the same time an analytic (smooth) linearizing coordinate change for $v$.

Now, assume that there exists an analytic (smooth) linearizing coordinate change for vector field $v$. Then, by Lemma \ref{triangular} there exists a linearizing coordinate change in the "triangular" form $\bar x = g(x, y), \bar y = y$. This implies, that $v(y) = \alpha y$ and $v(g) = g + \beta y$. Substituting this coordinate change into the \eqref{swipe} we get, that is coordinate change is, in fact, a linearizing coordinate change for Nijenhuis operator $R$. $\blacksquare$

Consider the following polynomial operator fields
\begin{equation*}
    \begin{aligned}
        \alpha = 0: & \quad R = \left(\begin{array}{cc}
y^2 & x \\
0 & y^2 \\
\end{array}\right), \quad F = \operatorname{det} R = y^4; \\
\alpha = r, 3 \leq r \in \mathbb N: & \quad R = \left(\begin{array}{cc}
     0 & x\\
     x^{\alpha - 1}& \alpha y \\ 
\end{array}\right), \quad F = \operatorname{det} R = - x^{\alpha}; \\
\alpha = - \frac{p}{q},\, p, q \in \mathbb N: & \quad R = \left(\begin{array}{cc}
0 & x + x^{p}y^{q} \\
0 & - \frac{p}{q}y \\
\end{array}\right); \\
\alpha = \frac{1}{m}, 2 \leq m \in \mathbb N: & \quad R = \left(\begin{array}{cc}
0 & x + \alpha y^{\frac{1}{\alpha}} \\
0 & \alpha y \\
\end{array}\right).
    \end{aligned}
\end{equation*}
The first operator is Nijenhuis by Corollary \ref{useful}, the third and the fourth are Nijenhuis by Corollary \ref{useful2}. The second operator field is Nijenhuis by Corollary \ref{local2} for $\operatorname{tr} R = \alpha y$ and $\operatorname{det} R = x^{\alpha}$. The linear part of the Taylor expansion of each operator field is $\mathfrak b_{1, \alpha}$. 

At the same time the first and the second operator field have no linearizing coordinate change in both smooth and analytic case by argument, similar to one in subsection \ref{romm}. The third and the fourth Nijenhuis operators does not have linearizing coordinate change by Lemma \ref{prenormal} as corresponding vector field $v$ from the condition of this Lemma by Theorem \ref{ilyash} can not be linearized (even in the formal category). This implies, obviously, that there are no analytic linearizing coordinate change. Theorem \ref{chen} implies, that there are no smooth linearizing coordinate change as well. Thus, we have shown, that $\mathfrak b_{1, \alpha}$ for $\alpha \in \Sigma_{\mathrm{an}}$ are degenerate in analytic category.

Now, assume, that $\alpha \notin \Sigma_{\mathrm{an}} \cup \Sigma_{\mathrm{u}}$. Consider two dimensional real plane with coordinates $x, y$ and Nijenhuis operator $R$, vanishing at the origin. We assume, that $\operatorname{tr} R = y$ and the linear part of the Taylor expansion of $R$ is
$$
\left(\begin{array}{cc}
     0 & \frac{1}{\alpha} x  \\
     0 & y 
\end{array}\right).
$$
This is the right-adjoint operator of $\mathfrak b_{1, \alpha}$, written in basis $\eta'_1 = \eta_1, \eta'_2 = \frac{1}{\alpha} \eta_2$ (here $\eta_1, \eta_2$ is basis from the Table 2 of Theorem \ref{class}). Consider $g = \operatorname{det} R$. We have that the Taylor expansion has no constant, linear and quadratic part. At the same time $\pd{R^1_2}{x} (0, 0) = \frac{1}{\alpha}$ and $\alpha \neq r, 3 \leq r \in \mathbb N$ and $r \neq - \frac{p}{q}$ with $p, q \in \mathbb N$. Thus, $R$ satisfies the conditions of Proposition \eqref{normal} and there exists a coordinate change, that transforms $R$ into
\begin{equation}\label{column}
R = \left(\begin{array}{cc}
     0 & f(x, y)  \\
     0 & y 
\end{array}\right).
    \end{equation}
Denote vector field $v = (f(x, y), y)^T$. Recall, that $\Omega$ stands for the set of negative Brjuno numbers. Assume, that $\alpha \neq 2$. Note, that together with the condition $\alpha \notin \Sigma_{\mathrm{an}} \cup \Sigma_{\mathrm{u}}$ this implies, that for $\alpha < 0, \alpha \in \Omega$ and for $\alpha > 0, \alpha \neq r, \frac{1}{r}$ with $2 \leq r \in \mathbb N$. By Corollary \ref{analyticc} there exists a linearizing coordinate change for vector field $v$. Thus, by Lemma \ref{prenormal} there exists a linearizing coordinate change for $R$.

Now, consider case $\alpha = 2$. By Theorem \ref{ilyash} there are three possible polynomial normal forms of $v$:
\begin{equation*}
    \begin{cases}
         \dot x = \frac{1}{2}x\\
         \dot y = y \\ 
    \end{cases} \, \, \, \, \, \, \quad
    \begin{cases}
         \dot x = \frac{1}{2}x\\
         \dot y = y - \frac{1}{2}x^2 \\ 
    \end{cases} \quad 
    \begin{cases}
         \dot x = x\\
         \dot y = 2y + \frac{1}{2}x^2 \\ 
    \end{cases}
\end{equation*}
By Theorem \ref{poincare} there exists an analytic coordinate change $\bar x = g(x, y), \bar y = h(x, y)$, that transforms $v$ into one of these forms. Note that in all three cases $v(g) = \frac{1}{2}g$. At the same time $v(y) = y$. 

Consider a pair of functions $g(x, y)$ and $y$. Denote the linearization operator of $v$ at the coordinate origin as $A$. Similar to the proof of Lemma \ref{triangular} we get, that the differentials $\ddd g (0, 0)$ and $\ddd y$ are eigenvectors for $A^*$ with eigenvalues $\frac{1}{2}$ and $1$ respectively. It is well-known from linear algebra, that these vectors are linearly independent. Thus, $\bar x = g(x, y), \bar y = y$ is a coordinate change and this is a linearizing coordinate change for $v$. By Lemma \ref{prenormal} there exists a linearizing coordinate change for $R$.

\subsection{Algebra $b_{1, \alpha}$ (smooth case)}

Consider linear vector field $v = (\alpha x, y)^T$ for $\alpha < 0$. Denote $s = - \frac{1}{\alpha}$ and denote $f(x, y)$ to be the smooth first integral of $v$ from Example \ref{smooth_int}. Define smooth function $g(x, y)$:
$$
g(x, y) = \begin{cases}
     \frac{x}{y} f(x, y) & xy \neq 0  \\
     0 & xy = 0.
\end{cases}
$$
For $xy \neq 0$ this function has the derivatives of all orders. Similar to the Example \ref{smooth_int} we define all the derivatives of $g$ to be zero on $xy = 0$. As a result we get a smooth function and by the definition of this function
\begin{equation}\label{idd}
g \pd{f}{x} = f\pd{f}{y}.     
\end{equation}
Consider operator field
$$
R = \left(\begin{array}{cc}
     f(x, y) & \alpha x + g(x, y)\\
     0 & y \\ 
\end{array}\right).
$$
We have $\operatorname{tr} R = f + y$, $\operatorname{det} R = yf$ and
\begin{equation*}
\begin{aligned}
       & \left( \begin{array}{cc}
             \pd{f}{x} ,&  \pd{f}{y} + 1  \\
        \end{array}\right)\left( \begin{array}{cc}
            y & - \alpha x - g \\
             0 & f 
        \end{array}\right) = \left( \begin{array}{cc}
             y\pd{f}{x} ,&  - \alpha x \pd{f}{x} - g \pd{f}{x} + f \pd{f}{y} + f \\
        \end{array}\right) = \\
        = &  \left( \begin{array}{cc}
             y\pd{f}{x} ,&  y\pd{f}{y} + f \\
        \end{array}\right) = \ddd \operatorname{det} R.
\end{aligned}
\end{equation*}
Here we used identity \eqref{idd} and $-\alpha x \pd{f}{x} = y \pd{f}{y}$, which follows from the fact, that $f$ is a first integral of $v = (\alpha x, y)^T$. By Proposition \ref{local1} this implies, that $R$ is Nijenhuis operator. At the same time for $\operatorname{det} R$ vanishes only on $xy = 0$ in the neighbourhood of the coordinate origin, while for righ-adjoint operator of $\mathfrak b_{1, \alpha}$ determinant vanishes identically. Thus, there are no linearizing coordinate change for $R$. For $\alpha = r, 3 \leq r \in \mathbb N$ and $\alpha = \frac{1}{r}, 2 \leq r \in \mathbb N$ the examples of Nijenhuis operators with no linearizing coordinate changes are identical to analytic case: one needs only to replace Theorem \ref{poincare} with Theorem \ref{chen} in the argument about the absence of linearizing coordinate changes. Thus, we have shown, that $\mathfrak b_{1, \alpha}$ for $\alpha \in \Sigma_{\mathrm{sm}}$ are degenerate in smooth category.

Same as in analytic case, consider two dimensional real plane with coordinates $x, y$ and smooth Nijenhuis operator $R$, vanishing at the origin. We assume, that $\operatorname{tr} R = y$ and the linear part of the Taylor expansion of $R$ is
$$
\left(\begin{array}{cc}
     0 & \frac{1}{\alpha} x  \\
     0 & y 
\end{array}\right).
$$
Consider $\alpha \notin \Sigma_{\rm{sm}}$. This implies, that $\alpha > 0$ and $\frac{1}{\alpha} \neq \frac{1}{r}$ with $3 \leq r \in \mathbb N$. Thus, by Proposition \ref{normal} there exists a smooth coordinate change, that transforms $R$ into the form \eqref{column}. 

Consider vector field $v = (f(x, y), y)^T$. Assume, that $\alpha \neq 2$. As $\alpha \notin \Sigma_{\mathrm{sm}}$, this implies, that $\alpha > 0, \alpha \neq r, \frac{1}{r}$ for $2 \leq r \in \mathbb N$. By Corollary \ref{smoothc} we get that there exists a smooth linearizing coordinate change for $v$. By Lemma \ref{prenormal} there exists a smooth linearizing coordinate change for $R$. The case $\alpha = 2$ is treated the same way as in analytic case with replacing Theorem \ref{poincare} with Theorem \ref{chen} in the argument.

\subsection{Algebra $\mathfrak b_3$}

Consider two dimensional real plane with coordinates $x, y$ and analytic (smooth) Nijenhuis operator $R$, vanishing at the origin. We assume, that $\operatorname{tr} R = y$ and the linear part of the Taylor expansion of $R$ is
$$
\left(\begin{array}{cc}
     0 & x + y\\
     0 & y 
\end{array}\right).
$$
We have, that $\pd{R^1_2}{x} (0, 0) = 1$ and the Taylor expansion of the determinant has no constant, linear and quadratic parts. Applying Proposition \ref{normal}, we get that there exists an analytic (smooth) coordinate change, that transforms $R$ into the form \ref{column}. 

Consider vector field $v = (f(x, y), y)^T$. The eigenvalues of linearization matrix are $\lambda_1 = \lambda_2 = 1$. Thus, by Corollaries \ref{analyticc} and \ref{smoothc} there exist an analytic (smooth) linearizing coordinate change. Thus, by Lemma \ref{prenormal} such coordinate change exists for operator field $R$.




\section{Appendix A: Vector fields on the plane}

\subsection{Linearization problem for vector fields}

Let $v$ be a vector field on two-dimensional real plane with coordinates $x, y$. Point ${\mathrm{p}}$ is \textbf{a critical point of $v$} if $v = 0$ at ${\mathrm{p}}$. In this section we always assume, that the coordinate origin is a critical point of $v$. 

Consider the Taylor expansion $v_1 + v_2 + \dots$ of vector field $v$ at the coordinate origin. The term $v_1$ is linear in coordinates, that is $v^1 = A^1_1 x + A^1_2 y, v^2 = A^2_1 x + A^2_2 y$. We call $A$ with components $A^i_j$ \textbf{the linearization operator} of $v$. It is a simple exercise to show, that $A^i_j$ define an operator on tangent space to the critical point.

Given vector field $v$ with non-zero linear part $v_1$ in its Taylor expansion one may ask if there a coordinate change around the critical point that transforms $v$ into linear vector field? This is the linearization problem for the vector fields. It was studies by many great mathematicians, starting with Poincare (see overview in \cite{ilya} and introduction in \cite{brjuno}). The answer to this question is given in terms of eigenvalues of $A^i_j$. 

In this section we assume, that eigenvalues of $A^i_j$ are real and non-zero $\lambda_1, \lambda_2$. The pair is called \textbf{resonant} if there exist $m, n \in \mathbb Z^+$, such that $m\lambda_1 + n \lambda_2 = \lambda_1$ and $m + n \geq 2$. Otherwise the pair is said to be \textbf{non-resonant}. With given restriction the pair $\lambda_1, \lambda_2$ is resonant if either $\frac{\lambda_1}{\lambda_2} = - \frac{p}{q}$, $\lambda_1 = r \lambda_2$ or $r \lambda_1 = \lambda_2$ with $p > q, r \geq 2 \in \mathbb N$. In particular $\lambda_1 = \lambda_2$ is non-resonant pair.

We distinguish linearization in three categories: formal, analytic and smooth. 

It is known, that (\cite{ilya}, Theorem 4) if the eigenvalues $\lambda_1, \lambda_2$ are non-resonant, then there exists a formal linearization. The proof of the theorem is based on a fact, that for non-resonant eigenvalues certain formal equation, known as homological equation, can be solved. Note, that this includes the case of $A^i_j$ being Jordan block with non-zero eigenvalues (see Lemma 4.5 \cite{ilya}, also \cite{basto}). In case of resonance one can obtain normal polynomial form.

\begin{theor}\label{ilyash}(\cite{ilya}, Proposition 4.29, Table 4.1)
Let $v$ be an analytic (or formal) vector field on the plane with critical point at the origin. Assume that eigenvalues of the linearization operator at this point are non-zero real $\lambda_1, \lambda_2$. Then there exists a formal coordinate change, that transforms $v$ into its polynomial normal form. The list of conditions and corresponding polynomial normal forms is given in the Table below:
\begin{equation*}
\begin{aligned}
& Table \, 5: \text{Polynomial normal forms of vector fields in dimension two} \\
    & \begin{tabular}{|c|c|c|}
    \hline
         Type & Conditions & Polynomial normal form  \\
         \hline
         No resonance & $\lambda_1, \lambda_2$ is non-resonant pair & linear \\
         \hline
         Resonant node I & $\lambda_1 = r \lambda_2$ and $r \in \mathbb N, r \geq 2$ & $\begin{aligned} &  \dot x = \lambda_1 x + a \lambda_2 y^r \\ & \dot y = \lambda_2 y, \\ & a = \{\pm 1, 0\}\end{aligned}$\\
         \hline
         Resonant node II & $r \lambda_1 = \lambda_2$ and $r \in \mathbb N, r \geq 2$ & $\begin{aligned} &  \dot x = \lambda_1 x \\ & \dot y = \lambda_2 y + a \lambda_1 x^r, \\ & a = \{\pm 1, 0\}\end{aligned}$\\
         \hline
         Resonant saddle & $\lambda_1 = - \frac{p}{q}\lambda_2, p, q \in \mathbb N$ & $\begin{aligned} & \dot x = \lambda_1 x \\ & \dot y = \lambda_2 y (1 \pm x^{qm}y^{pm} + b x^{2qm}y^{2pm}), \\ & b \in \mathbb R, m \in \mathbb N \end{aligned}$ \\
         \hline
    \end{tabular} \\
    \end{aligned}
\end{equation*}
The parameters $a, r$ for resonant node and $b, m$ for resonant saddle are formal invariants of normal forms.
\end{theor}

A critical point ${\mathrm{p}}$ of a vector field $v$ is called \textbf{elementary} if every eigenvalue of the linearization operator at ${\mathrm{p}}$ has a non-vanishing real part. The following Theorem deals with the linearization in smooth category. 
\begin{theor}\label{chen}
(K.T.Chen, \cite{chen}) Let $v$ and $w$ be two smooth vector fields having coordinate origin as an elementary critical point. Denote by $v_1 + v_2 + ...$ and $w_1 + w_2 + ...$ the respective Taylor's expansions of $v$ and $w$. Then there exists a $C^{\infty}$ transformation about $0$, which carries $v$ to $w$ if and only if there exists a formal transformation which carries the formal vector field $v_1 + v_2 + ...$ to $w_1 + w_2 + ...$.
\end{theor}

This theorem states, that linearization problem in smooth category for elementary vector fields is the same as formal linearization problem for their Taylor expansions. The Corollary holds.

\begin{corr}\label{smoothc}
Consider smooth vector field $v = (f(x, y), y)^T$ with critical point at the origin and eigenvalues of the linearization operator at this point being $\alpha, 1$. Assume, that 
$$
\alpha > 0, \quad \alpha \neq r, \frac{1}{r} \text{ with } 2 \leq r \in \mathbb N.
$$
Then there exists a smooth linearizing coordinate change for $v$.
\end{corr}
{\it Proof. } For a given $\alpha$ the polynomial normal form of $v$ in the Table 5 of Theorem \ref{ilyash} is linear. From Theorem \ref{chen} it follows, that there exists a smooth linearizing coordinate change. $\blacksquare$

In analytic category the linearization problem for vector fields is more complicated. The following theorem is a corollary of the classical Poincare-Dulac theorem. 

\begin{theor}\label{poincare}(\cite{ilya}, follows from Theorem 5.5)
Consider vector field $v$ on the plane with critical point at the origin and linearization operator with real positive non-zero eigenvalues. Then there exists an analytic coordinate change, that transforms $v$ into its polynomial normal form from Theorem \ref{ilyash}. 
\end{theor}

Recall, that $\Omega$ denotes the set of negative Brjuno numbers. We get the following theorem:

\begin{theor}\label{brjn} (\cite{brjuno}, follows form Theorem II, see also \cite{ilya} Theorem 5.22)
Consider analytic vector field $v$ on the plane with critical point at the origin and linearization operator with real non-zero eigenvalues $\lambda_1, \lambda_2$. Assume, that $\frac{\lambda_1}{\lambda_2} \in \Omega$. Then there exists a linearizing coordinate change.
\end{theor}

The following corollary holds.

\begin{corr}\label{analyticc}
Consider analytic vector field $v = (f(x, y), y)^T$ with critical point at the origin and eigenvalues of the linearization operator at this point being $\alpha, 1$. Assume, that either $\alpha \in \Omega$ or $\alpha > 0$ and $\alpha \neq r, \frac{1}{r}$ with $2 \leq r \in \mathbb N$. Then there exists an analytic linearizing coordinate change for $v$.
\end{corr}
{\it Proof. }If $\alpha \in \Omega$, then the analytic linearizing coordinate change exists by Theorem \ref{brjn}. If $\alpha > 0$ and $\alpha \neq r, \frac{1}{r}$ with $2 \leq r \in \mathbb N$, then the polynomial normal form from Theorem \ref{ilyash} is linear. By Theorem \ref{poincare} there exists a linearizing coordinat change. $\blacksquare$

In the end of this subsection we prove the following Lemma.

\begin{lemma}\label{triangular}
Consider $\alpha \neq 0$ and analytic (smooth) vector field $v = (f(x, y), y)$ on the plane with coordinates $x, y$ with critical point at the origin. Assume, that eigenvalues of linearization operator are $\alpha, 1$ and there exists a linearizing coordinate change around this point. Then there exists a linearizing coordinate change in "triangular" form $\bar{x} = g(x, y), \bar{y} = y$.
\end{lemma}
{\it Proof. } Vector field $v$ is an operator, that acts on analytic (smooth) functions by differentiation. Conditions of the Lemma imply that $v(y) = y$. The idea of the proof is to find among the functions a good pair for $y$ to define a proper coordinate change. We have three cases.

\textbf{Case $\alpha \neq 1$. } In this case the linearization operator $A$ can be brought to diagonal form with $\alpha, 1$ on diagonal. Consider the linearizing coordinate change, that exists by the condition of the Lemma. Assume that in this change $\bar{x} = g(x, y)$. We have that $v(g) = \alpha g$ and $\ddd g \neq 0$ at the origin. We have
$$
\alpha (\ddd g)_j = (\ddd v(g))_j = \pd{v^\alpha}{x^j} \pd{g}{x^{\alpha}} + v^{\alpha} \frac{\partial^2 g}{\partial x^{\alpha} x^j}.
$$
Substituting $x = y = 0$ into the equation we get, that $\ddd g$ at the origin is an eigenvector of $A^*$ with eigenvalue $\alpha$. From linear algebra we know, that eigenvectors for different eigenvalues are linearly independent. Thus, $\ddd g$ and $\ddd y$ are linearly independent at the origin and define locally analytic (smooth) coordinate change $\bar{x} = g(x, y), \bar{y} = y$. As $v(g) = \alpha g, v(y) = y$, we get that this is a linearizing coordinate change.

\textbf{Case $\alpha = 1$ and $A = \operatorname{Id}$}. Consider the analytic (smooth) linearizing coordinate change $\bar{x} = h_1(x, y), \bar{y} = h_2 (x, y)$ that exists by the condition of the Lemma. We have that $\ddd h_1$ and $\ddd h_2$ are linearly independent at the origin. Thus, at least one of these differentials is linearly independent with $\ddd y$. W.l.o.g. assume that this is $h_1$. We define analytic (smooth) coordinate change $\bar{x} = g(x, y) = h_1(x, y), \bar{y} = y$. By definition $v(g) = g, v(y) = y$, thus, this is a linearizing coordinate change.

\textbf{Case $\alpha = 1$ and $A$ is a Jordan block}. Consider the linearizing coordinate change, that exists by the conditions of the Lemma. In new coordinates $\bar{x}, \bar{y}$ we have $v = (\bar x + \bar y, \bar y)^T$. 

Consider analytic (smooth) solution $h$ of the equation $v(h) = h$ with extra condition $h(0, 0) = 0$. In local coordinates equation is 
$$
(\bar x + \bar y) \pd{h}{\bar x} + \bar y \pd{h}{\bar y} = h.
$$
Differentiating the equation in $\bar x$ we get
$$
0 = (\bar x + \bar y) \frac{\partial^2 h}{\partial \bar x \partial \bar x} + \bar y \frac{\partial^2 h}{\partial \bar y \partial \bar x} = v \Big( \pd{h}{\bar x}\Big).
$$
By \ref{integral2} any analytic (smooth) integral of $v = (\bar x + \bar y, \bar y)^T$ is constant around the origin. Thus, $v\big(\pd{h}{\bar x}\big) = 0$ implies, that $h = \gamma \bar x + k(\bar y)$ for some analytic (smooth) function of one variable $k$ and constant $\gamma$. 

For such $h$ equation $h = v(h)$ takes form $\gamma \bar x + k(\bar y) = \gamma (\bar x + \bar y) + \bar y k'(\bar y)$. Differentiating the latter in $\bar y$ yields
$$
\gamma + \bar y k''(\bar y) = 0.
$$
Substituting $\bar y = 0$ we get that $\gamma = 0$. This implies $\bar y k'' \equiv 0$ around the origin, thus, $k''(\bar y)$ identically vanish for $\bar y \neq 0$. By continuity it vanishes on the entire neighbourhood of the coordinate origin. As $h(0,0) = 0$ then $h = \gamma \bar y$ for some constant $\gamma \neq 0$. 

We have shown, that the space of solutions of equation $v(h) = h$ around the origin with condition $h(0, 0) = 0$ is one-dimensional. At the same time $v(y) = y$. Thus, the linearizing coordinate change, that exists by the conditions of Lemma has the form $\bar x = \bar{g}(x, y), \bar{y} = \beta y$ for some $\beta \neq 0$. Taking $\bar x = \frac{1}{\beta} \bar{g} (x, y) = g(x, y), \bar y = y$ yields the triangular linearizing analytic (smooth) coordinate change. The Lemma is proved. $\blacksquare$

\subsection{First integrals around critical points}

Assume that $v$ is vector field on the plane with critical point at the origin and linear part $v_1$ with linearization operator $A^i_j$. Recall, that the first integral of a vector field $v$ is function, such that $v(f) = 0$. We start with the following example.

\begin{example}\label{smooth_int}
Consider linear vector field $v = (\alpha x, y)$ for $\alpha < 0$. Fix the constant $s = - \frac{1}{\alpha} \in \mathbb R^+$. Define the function $f$ as follows:
$$
f(x, y) = \begin{cases}
     \exp \big(- \frac{1}{x^{2s} y^{2}}\big), & \quad xy \neq 0  \\
     0, & \quad xy = 0.
\end{cases}
$$
Defining all the partial derivatives of $f(x, y)$ as zero on the coordinate cross $xy = 0$ we obtain a function that is smooth on the entire plane and flat on the $xy = 0$. The partial derivatives of $h(x, y)$ satisfy the following identities:
$$
\pd{f}{x}(x, y) = \frac{2s}{x^{2s + 1} y^{2}}f(x, y), \quad \pd{f}{y} (x, y) = \frac{2}{x^{2s} y^{3}} f(x, y).
$$
In particular 
$$
\alpha x\pd{f}{x} + y \pd{f}{y} = 0 
$$
and, thus, for $\alpha < 0$ function $f(x, y)$ defines a smooth integral for linear vector field $v = (x, \alpha y)$.
\end{example}

The following two lemmas deal with integrals of the linear vector fields around the origin for dimension two.

\begin{lemma}\label{integral2}
Consider a real plane with coordinates $x, y$ and linear vector field $v = (x + y, y)^T$. It has no non-constant integrals in a neighbourhood of the origin in both smooth and analytic case.
\end{lemma}
{\it Proof. } The ODE $\dot x = x + y, \dot y = y$ can be explicitly integrated: 
\begin{equation*}
\begin{aligned}
    & x(t) = c_2 t \exp{t} +  c_1 \exp{t}, \\
    & y(t) = c_2 \exp{t},
\end{aligned}
\end{equation*}
were $c_1, c_2$ are arbitrary constants. 

We have $\lim_{t \to - \infty} x(t) = \lim_{t \to - \infty} y(t) = 0$. As a result, the closure of every integral curve contains the coordinate origin. The integral must be constant along the curves, thus, by continuity the value of integral on each curve coincides with value of the integral at the origin. As integral curves pass through every point of the neighouhood of the coordinate origin, we get that the first integral of $v$ is constant. $\blacksquare$

\begin{lemma}\label{integral1}
Consider a real plane with coordinates $x, y$ and linear vector field $v = (\alpha x, y)^T$ with $\alpha \neq 0$. Then
\begin{enumerate}
    \item It has non-constant smooth first integral on the neighbourhood of the origin if and only if $\alpha < 0$
    \item It has non-constant analytic first integral on the neighbourhood of the origin if and only if $\alpha = - \frac{p}{q}$ with $p, q \in \mathbb N$
\end{enumerate}
\end{lemma}
{\it Proof.} We start with smooth category. The ODE $\dot x = \alpha x, \dot y = y$ can be explicitly integrated: 
\begin{equation*}
    \begin{aligned}
        & x(t) = c_1 \exp{\alpha t}, \\
        & y(t) = c_2 \exp{t},
    \end{aligned}
\end{equation*}
where $c_1, c_2$ are arbitrary constants. Similar to the proof of \ref{integral2} the integral curves for $\alpha > 0$ tend to the coordinate origin as $t$ approaches $- \infty$. Thus, again, the closure of each curve contains the coordinate origin and every smooth coordinate is constant. For $\alpha < 0$ one takes $f(x, y)$ from Example \ref{smooth_int} to be non-constant smooth first integral of $v$.

Now consider analytic category and decomposition of the first integral $f$ of $v$ into the converging series $f_k + f_{k + 1} + \dots$. Here $f_k$ is the first non-zero term. Note, that $v(f) = v(f_k) + v(f_{k + 1}) + \dots$, where $v(f_i)$ is homogeneous polynomial of degree $i$. Thus, $v(f) = 0$ if and only if $v(f_i) = 0$ for all $i \geq k$. 

Consider $V^i$ to be the linear space of monomials $x^p y^q$ of degree $i = p + q$. We have
$$
v(x^p y^q) = (p\alpha + q) x^p y^q.
$$
This implies, that in the standard monomial basis $v$ acts on $V^i$ as a diagonal operator with eigenvalues $(p\alpha + q)$. The first integrals lie in the kernel of such operator. The kernel is non-empty if and only if at least one of eigenvalues is zero, that is $p + \alpha q = 0$ or $\alpha = - \frac{p}{q}$ for some $p, q \in \mathbb N$. 

Thus, if $\alpha \neq - \frac{p}{q}$ there are no non-constant first integral. If $\alpha = - \frac{p}{q}$ one takes $f = x^p y^q$ to be non-constant integral of the corresponding vector field $v$. The Lemma is proved. $\blacksquare$

The following Corollary holds.

\begin{corr}\label{nodes}
(Necessary condition for existing of the first integral) On the plane with coordinates $x, y$ consider analytic (smooth) vector field $v$ with critical point at the origin. Assume, that the eigenvalues of linearization operator at this point are non-zero real $\lambda_1, \lambda_2$. Assume, also, that there exists a non-constant analytic (smooth) first integral of $v$ around the origin. Then
\begin{enumerate}
    \item In analytic category $\frac{\lambda_1}{\lambda_2} = - \frac{p}{q}$ for $p \geq q \in \mathbb N$,
    \item In smooth category $\frac{\lambda_1}{\lambda_2} < 0$.
\end{enumerate}
\end{corr}
{\it Proof. } In analytic category consider the decompositions $v = v_1 + v_2 + \dots$ and $f = f_k + f_{k + 1} + \dots$ for vector field $v$ and first integral $f$. Here $f_k, k \geq 1$ is first non-zero term in the Taylor series of $f$. The formula for $v(f)$ takes form
$$
v(f) = v_1(f_k) + \big( v_2(f_k) + v_1(f_{k + 1})\big) + \dots
$$
On r.h.s. of this equation the term of degree $k$ is $v_1(f_k)$. As $v(f) = 0$ we get, that $f_k$ is the first integral of $v_1$. 

If $\lambda_1 \neq \lambda_2$, then w.l.o.g. the linearization operator $A$ can be considered diagonal. Taking $\frac{1}{\lambda_1} v_1$ instead of $v_1$ (rescaling by constant does not affect the fact that $f_k$ is the first integral of $v_1$) we get vector field in the form $v_1 = (\alpha x, y)^T$ with $\alpha = \frac{\lambda_1}{\lambda_2}$. Applying Lemma \ref{integral1} we get, that $\frac{\lambda_1}{\lambda_2} = - \frac{p}{q}$.

Now assume, that $\lambda_1 = \lambda_2$. W.l.o.g. we get that $v_1$ is either $v_1 = (x , y)^T$ or $(x + y, y)^T$. Applying Lemmas \ref{integral1} and \ref{integral2} we get the statement of the Corollary in analytic category.

Now, consider smooth category. Assume, that $\frac{\lambda_1}{\lambda_2} > 0$ and $\lambda_1, \lambda_2$ is non-resonant pair. By Theorem \ref{ilyash} the polynomial normal form of $v$ is linear. By Theorem \ref{chen} there exists a coordinate change, that transforms $v$ into linear vector field. W.l.o.g. we get that $\frac{1}{\lambda_1} v$ is either $(\alpha x, y)^T$ or $(x + y, y)$ depending on $\frac{\lambda_1}{\lambda_2}$. In both cases by Lemmas \ref{integral1} and \ref{integral2} there are no first integrals.

Now assume, that $\lambda_1, \lambda_2$ is a resonant pair of the first type $\lambda_1 = r \lambda_2$. Note, that any first integral of $v$ is a first integral of $\frac{1}{\lambda_2} v$. By Theorem \ref{ilyash} and Theorem \ref{chen} there exists a smooth coordinate change, that transforms $\frac{1}{\lambda_2} v$ into $v = (rx + a y^r, y)$. Here $r = \frac{\lambda_1}{\lambda_2}$ and $a$ is either $0, 1$ or $- 1$. The system of ODE $\dot x = rx + a y^r, \dot y = y$ can be explicitly integrated:
\begin{equation*}
    \begin{aligned}
        & x(t) = a c_2 t \exp rt + c_1 \exp rt, \\
        & y(t) = c_2 \exp t.
    \end{aligned}
\end{equation*}
We get that when $\lim_{t \to - \infty} x(t) = \lim_{t \to -\infty} y(t) = 0$. Following the same argument as in proof of Lemma \eqref{integral2}, we get, that the first integral in this case is constant around the origin. The case of the resonant node of the second type, that is $r \lambda_1 = \lambda_2$ is treated the same way. $\blacksquare$ 

\section{Appendix B: Morse lemma depending on parameters}

The following Lemma is well-known (see \cite{horm}, p.502). To keep our work self-sufficient, we provide it with proof. Recall, that the point ${\mathrm{p}}$ is said to be \textbf{critical} for a function $f(x, y)$ on the plane with coordinates $x, y$ if $\pd{f}{x} \vert_{\mathrm{p}}  = \pd{f}{y} \vert_{\mathrm{p}} = 0$.

\begin{lemma}\label{morse}
(Parametric Morse Lemma) Consider an analytic (smooth) function $f(x, y)$ on real plane with critical point at the coordinate origin. Assume, that
$$
\frac{\partial^2 f}{\partial x^2} (0, 0) = \gamma \neq 0, \quad \frac{\partial^2 f}{\partial x \partial y} (0, 0) = 0.
$$
There exists an analytic (smooth) coordinate change in the form $\bar{x} = g(x, y), \bar{y} = y$ with $g(0, 0) = 0$, such that
\begin{equation}\label{morse_normal}
f(\bar x, \bar y) = \operatorname{sgn}(\gamma) \, \bar x^2 + k(\bar y),    
\end{equation}
where $h$ is analytic(smooth) function of one variable with 
$$
k(0) = k'(0) = 0, \quad k''(0) = \frac{\partial^2 f}{\partial y^2} (0, 0).
$$
\end{lemma}
{\it Proof. } As $\frac{\partial^2 f}{\partial x^2} (0, 0) \neq 0$, then by the Implicit Function Theorem there exists an analytic (smooth) curve $r(y)$ such that 
$$
\pd{f}{x} (r(y), y) \equiv 0, \quad r (0) = 0.
$$
Function $f(x, y)$ can be written in the following integral form
\begin{equation*}
\begin{aligned}
    f(x, y) & = f(r(y), y) + \int \limits_{0}^1 \frac{\ddd}{\ddd t} f \big( tx + (1 - t) r(y), y \big) \ddd t = \\
    & = f(r(y), y) + \big( x - r(y)\big) \int \limits_{0}^1 \pd{f}{x} \big( tx + (1 - t) r(y), y \big) \ddd t
\end{aligned}
\end{equation*}
We denote
$$
\Phi(x, y) = \int \limits_{0}^1 \pd{f}{x} \big( tx + (1 - t) r(y), y \big) \ddd t
$$
By definition of $r(y)$ we have $\Phi(r(y), y) =\pd{f}{x}(r(y), y) \equiv 0$. For $\Phi$ the following integral formula formula holds
$$
\Phi(x, y) = \Phi(r(y), y) + \int \limits_{0}^1 \frac{\ddd}{\ddd t}\Phi \big( tx + (1 - t)r(y), y\big) = \big(x - r(y)\big) \int\limits_{0}^1 \pd{\Phi}{x} \big( tx + (1 - t)r(y), y\big)
$$
Denoting $F(x, y) = \int\limits_{0}^1 \pd{\Phi}{x} \big( tx + (1 - t)r(y), y\big)$ and substituting the expression into the original integral formula for $f$ we get
$$
f(x, y) = f(r(y), y) + \big( x - r(y)\big)^2 F(x, y).
$$
In this formula $\frac{\partial^2 f}{\partial x^2} (0, 0) = 2 F(0, 0) = \gamma \neq 0$. Consider functions
\begin{equation}\label{change}
    \bar x = (x - r(y)) \sqrt{|F(x, y)|} = g(x, y), \quad \bar y = y.
\end{equation}
These are analytic (smooth) functions in a neighbourhood of coordinate origin and $g(0, 0) = 0$. The Jacobian of $\bar x, \bar y$ at the origin is $\sqrt{\frac{|\gamma|}{2}} \neq 0$, thus, \eqref{change} defines a coordinate change in a neighbourhood of coordinate origin. In new coordinates after renaming $f(r(\bar y), \bar y) = k(\bar y)$ function $f$ takes the form \eqref{morse_normal}. By definition of curve $r(y)$ we get
$$
k'(y) = \pd{f}{x} (r(y), y) r'(y) + \pd{f}{y}(r(y), y) = \pd{f}{y}(r(y), y), \quad k''(y) = \frac{\partial^2 f}{\partial x \partial y} r'(y) + \frac{\partial^2 f}{\partial y^2}.
$$
This implies, that $k(0) = k'(0) = 0$ and $k''(0) = \frac{\partial^2 f}{\partial y^2} (0, 0)$. The Lemma is proved. $\blacksquare$.

\begin{corr}\label{morse2}
Consider an analytic (smooth) function $f(x, y)$ on real plane with critical point at the coordinate origin. Assume, that $f(0, 0) = 0$ and
$$
\frac{\partial^2 f}{\partial x^2} (0, 0) = \frac{\partial^2 f}{\partial y^2} (0, 0) = \frac{\partial^2 f}{\partial x \partial y} (0, 0) = 0.
$$
In other words Taylor expansion of $f$ at the coordinate origin has no constant, linear and quadratic part. For every $\alpha \neq 0$ there exist analytic (smooth) functions $h(x, y)$ and $k(y)$ such, that 
$$
h(0, 0) = 0, \quad \Big(\pd{h}{y}(0, 0)\Big)^2 = \frac{\alpha^2}{4} \neq 0, \quad \pd{h}{x} (0, 0) = 0, \quad k(0) = k'(0) = k''(0) = 0
$$
and
\begin{equation}\label{morse2_normal}
f(x, y) = \frac{\alpha^2}{4} y^2 - h^2(x, y) + k(x).
\end{equation}
\end{corr}
{\it Proof. } Consider function $\bar f = f - \frac{\alpha^2}{4} y^2$. For $\bar f$ the coordinate origin is still a critical point. At the same time $\frac{\partial^2 \bar f}{\partial y^2} (0, 0) = - \frac{\alpha^2}{2} \neq 0$ and $\frac{\partial^2 f}{\partial x \partial y} (0, 0) = 0$.

Applying Lemma \ref{morse} (in our situation the coordinates $x, y$ are interchanged with respect to the coordinates in the statement of the Lemma) we get that there exists a coordinate change $\bar x = x, \bar y = g(x, y)$, such that $\bar f = - \bar y^2 + k(\bar x)$. We also have, that $g(0, 0) = 0$. In old coordinates this formula takes the form
\begin{equation}\label{tr1}
f - \frac{\alpha^2}{4} y^2 = - g^2(x, y) + k(x).    
\end{equation}
Renaming $g(x, y)$ as $h(x, y)$ we get the formula \eqref{morse2_normal} with condition $h(0, 0) = 0$. Differentiating both sides of the equation \eqref{tr1} we get
\begin{equation*}
    \begin{aligned}
        & \pd{f}{x} = - 2 h \pd{h}{x} + k'(x), \\
        & \frac{\partial^2 f}{\partial x^2} = - 2 \Big( \pd{h}{x}\Big)^2 - 2 h \frac{\partial^2 h}{\partial x^2} +  k''(x), \\
        & \frac{\partial^2 f}{\partial x \partial y} = - 2 \pd{h}{x}\pd{h}{y} - 2 h \frac{\partial^2 h}{\partial x \partial y}, \\
        & \frac{\partial^2 f}{\partial y^2} - \frac{\alpha^2}{2} = - 2 \Big( \pd{h}{y}\Big)^2 - 2 h \, \frac{\partial^2 h}{\partial y^2}. 
    \end{aligned}
\end{equation*}
Substituting $x = y = 0$ we get $\big(\pd{h}{y}(0, 0)\big)^2 = \frac{\alpha^2}{4} \neq 0$ from the fourth equation. After that the third equation yields $\pd{h}{x}(0, 0) = 0$. Substituting both into the first and second one, we get $k'(0) = k''(0) = 0$. The Corollary is proved. $\blacksquare$


\end{document}